\documentclass{amsproc}%
\usepackage{amsfonts}
\usepackage{amsmath}
\usepackage{amssymb}
\usepackage{graphicx}%
\theoremstyle{plain}

\newtheorem{corollary}{Corollary}

\newtheorem{definition}{Definition}

\newtheorem{lemma}{Lemma}

\newtheorem{proposition}{Proposition}

\newtheorem{theorem}{Theorem}
\numberwithin{equation}{section}
\begin{document}
\title[Type $B_{2}$ Intertwining]{An Intertwining Operator for the Group $B_{2}$}
\author{Charles F. Dunkl}
\address{Department of Mathematics, PO Box 400137, University of Virginia,
Charlottesville, VA 22904-4137}
\email{cfd5z@virginia.edu}
\urladdr{http://www.people.virginia.edu/\symbol{126}cfd5z/}
\subjclass[2000]{Primary 33C80, 33C20; Secondary 33C70, 43A80}
\keywords{intertwining operator, Dunkl operators, balanced hypergeometric functions}
\dedicatory{ }
\begin{abstract}
There is a commutative algebra of differential-difference operators, acting on
polynomials on $\mathbb{R}^{2}$, associated with the reflection group $B_{2}$.
This paper presents an integral transform which intertwines this algebra,
allowing one free parameter, with the algebra of partial derivatives. The
method of proof depends on properties of a certain class of balanced
terminating hypergeometric series of $_{4}F_{3}$-type. These properties are in
the form of recurrence and contiguity relations and are proved herein.

\end{abstract}
\date{July 31, 2006}
\maketitle

\section{Introduction}

\subsection{Overview}

We construct an integral for the intertwining operator $V$ associated to the
reflection group of type $B_{2}$ (order 8) acting on $\mathbb{R}^{2}$, with
\textit{one }parameter $\kappa$. For polynomials or adequately smooth
functions in $x=\left(  x_{1},x_{2}\right)  $ define the
differential-difference operators:
\begin{align}
T_{1}f\left(  x\right)   &  =\frac{\partial}{\partial x_{1}}f\left(  x\right)
+\kappa_{1}\frac{f\left(  x\right)  -f\left(  -x_{1},x_{2}\right)  }{x_{1}%
}\label{ddops}\\
&  +\kappa_{2}\frac{f\left(  x\right)  -f\left(  x_{2},x_{1}\right)  }%
{x_{1}-x_{2}}+\kappa_{2}\frac{f\left(  x\right)  -f\left(  -x_{2}%
,-x_{1}\right)  }{x_{1}+x_{2}},\nonumber\\
T_{2}f\left(  x\right)   &  =\frac{\partial}{\partial x_{1}}f\left(  x\right)
+\kappa_{1}\frac{f\left(  x\right)  -f\left(  x_{1},-x_{2}\right)  }{x_{2}%
}\nonumber\\
&  +\kappa_{2}\frac{f\left(  x\right)  -f\left(  x_{2},x_{1}\right)  }%
{x_{2}-x_{1}}+\kappa_{2}\frac{f\left(  x\right)  -f\left(  -x_{2}%
,-x_{1}\right)  }{x_{1}+x_{2}}.\nonumber
\end{align}
These operators are special cases of those defined by the author in \cite{D1}.
Their key property is commutativity, $T_{1}T_{2}=T_{2}T_{1}$. We deal only
with the restricted case $\kappa_{1}=\kappa_{2}=\kappa$. The intertwining
operator $V$ preserves the degree of homogeneous polynomials and satisfies
$V\left(  \frac{\partial}{\partial x_{i}}\right)  f\left(  x\right)
=T_{i}Vf\left(  x\right)  $ for $i=1,2$, and $V1=1$. In \cite{DJO} it was
shown that $V$ exists and is one-to-one as a map on polynomials for any
$\kappa$ except for the set $\left\{  -\frac{m}{4}:m\in\mathbb{N},\frac{m}%
{4}\notin\mathbb{Z}\right\}  $ of singular values. The definition and
existence of $V$ for $\kappa>0$ was shown in \cite{D2}. Later R\"{o}sler
\cite{R1} (see also \cite{R2}) proved that the functional $f\mapsto Vf\left(
x\right)  $ is given by integration with respect to a positive measure, for
each $x\in\mathbb{R}$. An explicit integral was found for the group $S_{3}$
(symmetric group on 3 objects) in \cite{D3}. Xu \cite{X} found an intertwining
transform for $B_{2}$ under restrictive conditions on degree and $\kappa
_{1},\kappa_{2}$. We will use a criterion, valid for any finite reflection
group, for an operator to be equal to $V$ which was proven in that paper.
There is a kernel which contains all the information about the action of $V$
on polynomials. For $x,y\in\mathbb{R}^{2}$ let $\left\langle x,y\right\rangle
:=\sum_{i=1}^{2}x_{i}y_{i}$, then
\[
K\left(  x,y\right)  =V^{x}\left(  \exp\left\langle x,y\right\rangle \right)
\]
is entire in $x$ and $y$ ($V^{x}$ acts on the variable $x$). Further let
$K_{n}\left(  x,y\right)  =\frac{1}{n!}V^{x}\left(  \left\langle
x,y\right\rangle ^{n}\right)  $ for $n\in\mathbb{N}_{0}:=\left\{
0,1,2,\ldots\right\}  $. The book of Dunkl and Xu \cite{DX} is a reference for
differential-difference operators, the intertwining operator and the kernel
$K\left(  x,y\right)  $ for any finite reflection group. Asymptotic formulae
for the kernel are discussed in \cite{R2}. The starting point for the $S_{3}$
result was an integral formula of Harish-Chandra involving integration over a
compact Lie group. This gives a special case ($\kappa=1$ and symmetrized) of
the intertwining operator for the associated Weyl group.

\subsection{The Symplectic Group Integral}

We use an approach similar to that in \cite[sect.2]{D3}. Consider $B_{2}$ as
the Weyl group of the (compact) symplectic group $Sp\left(  2\right)  $. This
group can be described as the group of $2\times2$ unitary block matrices%
\[
U=\left(
\begin{array}
[c]{cc}%
A & B\\
-\overline{B} & \overline{A}%
\end{array}
\right)  ,UU^{\ast}=I,
\]
where $A,B$ are $2\times2$ complex matrices. The subgroup%
\[
T=\left\{  \mathrm{diag}\left(  e^{\mathrm{i}\theta_{1}},e^{\mathrm{i}%
\theta_{2}},e^{-\mathrm{i}\theta_{1}},e^{-\mathrm{i}\theta_{2}}\right)
:\theta_{1},\theta_{2}\in\mathbb{R}\right\}
\]
is a maximal torus (\textquotedblleft diag\textquotedblright\ denotes the
$4\times4$ diagonal matrix with the specified entries). Identify
$\mathbb{R}^{2}$ with the complexification of the Lie algebra of $T$ by the
map $\delta:x\mapsto2^{-1/2}\mathrm{diag}\left(  x_{1},x_{2},-x_{1}%
,-x_{2}\right)  $ (the purpose of the factor $2^{-1/2}$ is to get $Tr\left(
\delta\left(  x\right)  \delta\left(  y\right)  \right)  =\left\langle
x,y\right\rangle $). The formula of Harish-Chandra (Helgason \cite[p.328]{H})
specializes to%
\[
\int_{Sp\left(  2\right)  }\exp\left(  Tr\left(  U\delta\left(  x\right)
U^{\ast}\delta\left(  y\right)  \right)  \right)  dm\left(  U\right)
=\frac{3}{2}\frac{\sum_{w\in B_{2}}\mathrm{\det}\left(  w\right)  \exp\left(
\left\langle xw,y\right\rangle \right)  }{p\left(  x\right)  p\left(
y\right)  },
\]
where $dm\left(  U\right)  $ is normalized Haar measure on $Sp\left(
2\right)  $ and $p\left(  x\right)  =x_{1}x_{2}\left(  x_{1}^{2}-x_{2}%
^{2}\right)  $. The right side of the formula is an expression for the kernel
$K^{0}\left(  x,y\right)  $ at $\kappa=1$. Thus the left side suggests a
construction of an integral formula for the intertwining operator. We compute
\[
Tr\left(  U\delta\left(  x\right)  U^{\ast}\delta\left(  y\right)  \right)
=\sum_{i=1}^{2}\sum_{j=1}^{2}x_{i}y_{j}\left(  \left\vert A_{ji}\right\vert
^{2}-\left\vert B_{ji}\right\vert ^{2}\right)  .
\]
Thus one needs to integrate functions of the four variables $\left(
\left\vert A_{ji}\right\vert ^{2}-\left\vert B_{ji}\right\vert ^{2}\right)  $
with respect to Haar measure. One applies integration over the subgroup
$Sp\left(  1\right)  $ and the homogeneous space $Sp\left(  2\right)
/Sp\left(  1\right)  $. Since this only gives the $\kappa=1$ situation more
experimentation is needed to make a conjecture about arbitrary $\kappa$. By
direct polynomial calculation we find the symmetrized kernel $K_{n}^{0}\left(
x,y\right)  =\frac{1}{8}\sum_{w\in B_{2}}V^{x}\left\langle xw,y\right\rangle
^{n}/n!$ for small $n$ ($\leq8$) and try powers of the trigonometric functions
to produce $K_{n}^{0}\left(  x,y\right)  =\frac{1}{n!}\int_{\Omega
}\left\langle x\tau\left(  q\right)  ,y\right\rangle ^{n}d\mu\left(  q\right)
$. This approach, however, has as yet not produced a solution for the
two-parameter situation.

\subsection{Group actions and the measure}

The group $B_{2}$ is generated by the reflections $\sigma_{1}:=\left(
\begin{array}
[c]{cc}%
-1 & 0\\
0 & 1
\end{array}
\right)  $ and $\sigma_{2}:=\left(
\begin{array}
[c]{cc}%
0 & 1\\
1 & 0
\end{array}
\right)  $. The intertwining operator has the form%
\begin{align*}
Vf\left(  x\right)   &  =\int_{\Omega}f\left(  x\tau\left(  q\right)  \right)
g\left(  q\right)  d\mu\left(  q\right)  ,\\
\tau\left(  q\right)   &  =\left(
\begin{array}
[c]{cc}%
q_{1} & q_{3}\\
q_{2} & q_{4}%
\end{array}
\right)  ,
\end{align*}
where $\mu\left(  q\right)  $ is two-sided invariant for $B_{2}$ and $g\left(
q\right)  $ is invariant under the action $w:\tau\left(  q\right)  \rightarrow
w\tau\left(  q\right)  w^{-1}$ for all $w\in B_{2}$ (actually $g\left(
q\right)  $ involves one more variable of integration). There are left and
right representations of $B_{2}$ on $q$ (that is,\ $\mathbb{R}^{4}$), defined
implicitly for $w\in B_{2}$ by%
\begin{align*}
\tau\left(  q\lambda\left(  w\right)  \right)   &  =w^{-1}\tau\left(
q\right)  ,\\
\tau\left(  q\rho\left(  w\right)  \right)   &  =\tau\left(  q\right)  w.
\end{align*}
For example $q\lambda\left(  \sigma_{1}\right)  =\left(  -q_{1},q_{2}%
,-q_{3},q_{4}\right)  $ and $q\rho\left(  \sigma_{2}\right)  =\left(
q_{3},q_{4},q_{1},q_{2}\right)  $. The invariance conditions are $d\mu\left(
q\lambda\left(  w\right)  \right)  =d\mu\left(  q\rho\left(  w\right)
\right)  =d\mu\left(  q\right)  $ and $g\left(  q\lambda\left(  w\right)
\rho\left(  w\right)  \right)  =g\left(  q\right)  $ for all $w\in B_{2}$. It
suffices to check invariance for the generating reflections, that is, $g$ and
$\mu$ must be invariant under $q\mapsto q\lambda\left(  \sigma_{1}\right)
\rho\left(  \sigma_{1}\right)  =\left(  q_{1},-q_{2},-q_{3},q_{4}\right)  $
and $q\mapsto q\lambda\left(  \sigma_{2}\right)  \rho\left(  \sigma
_{2}\right)  =\left(  q_{4},q_{3},q_{2},q_{1}\right)  $, and additionally
$\mu$ must be invariant under $q\mapsto q\lambda\left(  \sigma_{1}\right)
=\left(  -q_{1},q_{2},-q_{3},q_{4}\right)  $ and $q\mapsto q\lambda\left(
\sigma_{2}\right)  =\left(  q_{2},q_{1},q_{4},q_{3}\right)  $.

The measure $\mu$ is actually an integral over six variables: let%
\begin{align*}
q_{1}  &  =u\cos\psi_{1}\\
q_{2}  &  =\left(  1-u\right)  \cos\psi_{2}\\
q_{3}  &  =\left(  1-u\right)  \left(  \cos\psi_{2}\cos\theta+\sin\psi_{2}%
\sin\theta\cos\phi_{2}\right) \\
q_{4}  &  =u\left(  \cos\psi_{1}\cos\theta+\sin\psi_{1}\sin\theta\cos\phi
_{1}\right)  ;
\end{align*}
the region $\Omega$ of integration is $0\leq u\leq1,0\leq\theta,\phi_{i}%
,\psi_{i}\leq\pi~\left(  i=1,2\right)  $. For $\kappa>\frac{1}{2}$ the measure
is%
\[
d\mu\left(  q\right)  =c_{\kappa}\left(  u\left(  1-u\right)  \sin\psi_{1}%
\sin\psi_{2}\sin\theta\right)  ^{2\kappa-1}\left(  \sin\phi_{1}\sin\phi
_{2}\right)  ^{2\kappa-2}dud\psi_{1}d\psi_{2}d\theta d\phi_{1}d\phi_{2}%
\]
times the normalizing constant (product of gamma functions) so $\int_{\Omega
}d\mu=1$. We can now state the main result for homogeneous polynomials: let
\[
D_{0}=\left(  q_{1}+q_{4}\right)  \left(  \frac{\partial}{\partial q_{1}%
}+\frac{\partial}{\partial q_{4}}\right)  -\left(  q_{2}-q_{3}\right)  \left(
\frac{\partial}{\partial q_{2}}-\frac{\partial}{\partial q_{3}}\right)  ,
\]
and suppose $f\left(  x\right)  $ is a homogeneous polynomial in $x$ of degree
$n$.

\begin{enumerate}
\item If $n$ is odd, then%
\begin{equation}
Vf\left(  x\right)  =\int_{\Omega}2\left(  q_{1}+q_{4}\right)  f\left(
x\tau\left(  q\right)  \right)  d\mu\left(  q\right)  . \label{oddfV}%
\end{equation}

\item If $n$ is even, then%
\begin{equation}
Vf\left(  x\right)  =\int_{\Omega}\left\{  1+q_{1}q_{4}-q_{2}q_{3}+\frac
{1}{4\kappa+n}D_{0}\right\}  f\left(  \tau\left(  q\right)  x\right)
d\mu\left(  q\right)  . \label{evenfV}%
\end{equation}

\end{enumerate}

Further, the Bessel function $K^{0}\left(  x,y\right)  =\frac{1}{8}\sum_{w\in
B_{2}}K\left(  xw,y\right)  $ is given by the positive integral:%
\[
K^{0}\left(  x,y\right)  =\int_{\Omega}\exp\left(  \left\langle x\tau\left(
q\right)  ,y\right\rangle \right)  d\mu\left(  q\right)  .
\]

We will consider the relation to the representations of $B_{2}$ and the
techniques of proving the formulae in the following sections. The required
invariance properties of $\mu$ are made clear by a change-of-variables.

\begin{lemma}
\label{invmu}In terms of $q,\theta,u$ and the auxiliary variables%
\begin{align*}
z_{1}\left(  q_{1},q_{4},u,\theta\right)   &  =\frac{u^{2}\sin^{2}\theta
-q_{1}^{2}-q_{4}^{2}+2q_{1}q_{4}\cos\theta}{\sin^{2}\theta},\\
z_{2}\left(  q_{2},q_{3},u,\theta\right)   &  =\frac{\left(  1-u\right)
^{2}\sin^{2}\theta-q_{2}^{2}-q_{3}^{2}+2q_{2}q_{3}\cos\theta}{\sin^{2}\theta}%
\end{align*}
the measure is given by
\[
d\mu=c_{\kappa}\left(  z_{1}z_{2}\right)  ^{\kappa-3/2}\left(  \sin
\theta\right)  ^{2\kappa-3}dq_{1}dq_{2}dq_{3}dq_{4}dud\theta,
\]
and the region of integration $\Omega_{1}$ is implicitly defined by $z_{1}%
\geq0,z_{2}\geq0,0\leq\theta\leq\pi,0\leq u\leq1$.
\end{lemma}

\begin{proof}
In terms of $q_{1},q_{2},q_{3},q_{4},u,\theta$ the Jacobian is%
\[
J=\dfrac{\partial\left(  q_{1},q_{2},q_{3},q_{4},u,\theta\right)  }%
{\partial\left(  \phi_{1},\phi_{2},\psi_{1},\psi_{2},u,\theta\right)
}=\left(  u\left(  1-u\right)  \sin\psi_{1}\sin\psi_{2}\sin\theta\right)
^{2}\sin\phi_{1}\sin\phi_{2}.
\]
Then
\[
d\mu=c_{\kappa}\left(  u\left(  1-u\right)  \sin\psi_{1}\sin\psi_{2}\sin
\phi_{1}\sin\phi_{2}\sin\theta\right)  ^{2\kappa-3}Jdud\psi_{1}d\psi
_{2}d\theta d\phi_{1}d\phi_{2}.
\]
Observe $u\sin\psi_{1}\sin\theta\cos\phi_{1}=q_{4}-q_{1}\cos\theta$ and
$u^{2}\sin^{2}\psi_{1}=u^{2}-q_{1}^{2}$ thus%
\begin{align*}
\sin^{2}\phi_{1}  &  =1-\left(  \frac{q_{4}-q_{1}\cos\theta}{u\sin\psi_{1}%
\sin\theta}\right)  ^{2},\\
\left(  u\sin\psi_{1}\sin\phi_{1}\right)  ^{2}  &  =\left(  u^{2}-q_{1}%
^{2}\right)  -\frac{\left(  q_{4}-q_{1}\cos\theta\right)  ^{2}}{\sin^{2}%
\theta}\\
&  =z_{1}\left(  q_{1},q_{4},u,\theta\right)  ,
\end{align*}
and similarly%
\begin{align*}
\left(  \left(  1-u\right)  \sin\psi_{2}\sin\phi_{2}\right)  ^{2}  &  =\left(
\left(  1-u\right)  ^{2}-q_{2}^{2}\right)  -\frac{\left(  q_{3}-q_{2}%
\cos\theta\right)  ^{2}}{\sin^{2}\theta}\\
&  =z_{2}\left(  q_{2},q_{3},u,\theta\right)  .
\end{align*}
Thus $d\mu=c_{\kappa}\left(  z_{1}z_{2}\right)  ^{\kappa-3/2}\left(
\sin\theta\right)  ^{2\kappa-3}dq_{1}dq_{2}dq_{3}dq_{4}dud\theta$.
\end{proof}

The measure $\mu$ is invariant under the transpositions $\left(  q_{1}%
,q_{4}\right)  ,\left(  q_{2},q_{3}\right)  $, and the involutions $q\mapsto
q\lambda\left(  \sigma_{1}\right)  =\left(  q_{2},q_{1},q_{4},q_{3}\right)  $
and $q\mapsto q\lambda\left(  \sigma_{1}\right)  =\left(  -q_{1},-q_{2}%
,q_{3},q_{4}\right)  $ because of the equations%
\begin{align*}
z_{1}\left(  q_{2},q_{3},u,\theta\right)  z_{2}\left(  q_{1},q_{4}%
,u,\theta\right)   &  =z_{1}\left(  q_{1},q_{4},1-u,\theta\right)
z_{2}\left(  q_{2},q_{3},1-u,\theta\right)  ,\\
z_{1}\left(  -q_{1},q_{4},u,\theta\right)  z_{2}\left(  -q_{2},q_{3}%
,u,\theta\right)   &  =z_{1}\left(  q_{1},q_{4},u,\pi-\theta\right)
z_{2}\left(  q_{2},q_{3},u,\pi-\theta\right)  .
\end{align*}
Similarly there is invariance under the right action of $B_{2}$, that is,
$q\mapsto q\rho\left(  \sigma_{i}\right)  ,i=1,2$.

To integrate the typical monomial $q^{\alpha}:=q_{1}^{\alpha_{1}}q_{2}%
^{\alpha_{2}}q_{3}^{\alpha_{3}}q_{4}^{\alpha_{4}}$ where $\alpha\in
\mathbb{N}_{0}^{4}$ ($q^{\alpha}$ is of degree $\left\vert \alpha\right\vert
:=\sum_{i=1}^{4}\alpha_{i}$) set
\begin{align*}
b_{0}  &  =\left(  \alpha_{2}+\alpha_{3}\right)  /2,\\
b_{1}  &  =\left(  \alpha_{1}+\alpha_{4}\right)  /2,\\
b_{2}  &  =\left(  \alpha_{2}+\alpha_{4}\right)  /2,\\
b_{3}  &  =\left(  \alpha_{3}+\alpha_{4}\right)  /2,
\end{align*}
and let $s\left(  \alpha_{1},\alpha_{2},\alpha_{3},\alpha_{4}\right)
:=\int_{\Omega}q^{\alpha}d\mu\left(  q\right)  $. By a parity argument the
nonzero integrals occur only for integer values of the $b_{i}$.

\begin{proposition}
The normalizing constant is%
\[
c_{\kappa}=4^{\kappa-1}\frac{\left(  2\kappa-1\right)  ^{2}\Gamma\left(
2\kappa+\frac{1}{2}\right)  }{\pi^{5/2}\Gamma\left(  \kappa\right)  ^{2}},
\]
and for $\alpha\in\mathbb{N}_{0}^{4}$ if $\alpha_{1}\equiv\alpha_{2}%
\equiv\alpha_{3}\equiv\alpha_{4}\operatorname{mod}2$ (all even or all odd)
then%
\begin{subequations}
\begin{align}
s\left(  \alpha_{1},\alpha_{2},\alpha_{3},\alpha_{4}\right)   &
=\frac{\left(  2\kappa\right)  _{2b_{1}}\left(  2\kappa\right)  _{2b_{0}%
}\left(  \frac{1}{2}\right)  _{b_{1}}\left(  \frac{1}{2}\right)  _{b_{0}%
}\left(  \frac{1}{2}\right)  _{b_{3}}}{\left(  4\kappa\right)  _{2b_{1}%
+2b_{0}}\left(  \kappa+\frac{1}{2}\right)  _{b_{1}}\left(  \kappa+\frac{1}%
{2}\right)  _{b_{0}}\left(  \kappa+\frac{1}{2}\right)  _{b_{3}}}%
\label{dblsum}\\
&  \times\sum_{i=0}^{\left\lfloor \alpha_{4}/2\right\rfloor }\sum
_{j=0}^{\left\lfloor \alpha_{3}/2\right\rfloor }\frac{\left(  -\alpha
_{4}\right)  _{2i}\left(  -\alpha_{3}\right)  _{2j}\left(  \kappa\right)
_{i+j}}{i!j!\left(  \frac{1}{2}-b_{1}\right)  _{i}\left(  \frac{1}{2}%
-b_{0}\right)  _{j}\left(  \frac{1}{2}-b_{3}\right)  _{i+j}}2^{-2i-2j}%
,\nonumber
\end{align}
otherwise $s\left(  \alpha_{1},\alpha_{2},\alpha_{3},\alpha_{4}\right)  =0$.
\end{subequations}
\end{proposition}

\begin{proof}
Expand $q^{\alpha}$ in terms of $u,\theta,\psi_{1},\psi_{2},\phi_{1},\phi_{2}$
with the binomial theorem and collect terms. The result is%
\begin{align*}
s\left(  \alpha_{1},\alpha_{2},\alpha_{3},\alpha_{4}\right)   &  =c_{\kappa
}\sum_{j=0}^{\alpha_{3}}\sum_{i=0}^{\alpha_{4}}\binom{\alpha_{3}}{j}%
\binom{\alpha_{4}}{i}\int_{0}^{1}u^{\alpha_{1}+\alpha_{4}+2\kappa-1}\left(
1-u\right)  ^{\alpha_{2}+\alpha_{3}+2\kappa-1}du\\
&  \times\int_{0}^{\pi}\left(  \cos\psi_{1}\right)  ^{\alpha_{1}+\alpha_{4}%
-i}\left(  \sin\psi_{1}\right)  ^{i+2\kappa-1}d\psi_{2}\int_{0}^{\pi}\left(
\cos\phi_{1}\right)  ^{i}\left(  \sin\phi_{1}\right)  ^{2\kappa-2}d\phi_{2}\\
&  \times\int_{0}^{\pi}\left(  \cos\psi_{2}\right)  ^{\alpha_{2}+\alpha_{3}%
-j}\left(  \sin\psi_{2}\right)  ^{j+2\kappa-1}d\psi_{1}\int_{0}^{\pi}\left(
\cos\phi_{2}\right)  ^{j}\left(  \sin\phi_{2}\right)  ^{2\kappa-2}d\phi_{1}\\
&  \times\int_{0}^{\pi}\left(  \cos\theta\right)  ^{\alpha_{3}+\alpha_{4}%
-i-j}\left(  \sin\theta\right)  ^{i+j+2\kappa-1}d\theta.
\end{align*}
Recall $\int_{0}^{\pi}\cos^{n}\theta\sin^{\lambda}\theta d\theta$ equals zero
if $n$ is odd and equals $B\left(  \frac{n+1}{2},\frac{\lambda+1}{2}\right)  $
if $n$ is even, for $n=0,1,2,\ldots$and $\lambda>-1$. For the respective
integrals to be nonzero $i$ and $j$ must be even because of the $\phi_{1}$ and
$\phi_{2}$ integrals and hence $\alpha_{2}+\alpha_{3},\alpha_{1}+\alpha
_{4},\alpha_{3}+\alpha_{4}$ must be even (in the $\psi_{1},\psi_{2},\theta$
integrals). So replace $i,j$ by $2i,2j$ respectively and assume that the
entries of $\alpha$ are all even or all odd. Then%
\begin{align*}
s\left(  \alpha_{1},\alpha_{2},\alpha_{3},\alpha_{4}\right)   &  =c_{\kappa
}\sum_{j=0}^{\left\lfloor \alpha_{3}/2\right\rfloor }\sum_{i=0}^{\left\lfloor
\alpha_{4}/2\right\rfloor }\binom{\alpha_{3}}{2j}\binom{\alpha_{4}}%
{2i}B\left(  \alpha_{1}+\alpha_{4}+2\kappa,\alpha_{2}+\alpha_{3}%
+2\kappa\right) \\
&  \times B\left(  \frac{\alpha_{1}+\alpha_{4}+1}{2}-i,i+\kappa\right)
B\left(  i+\frac{1}{2},\kappa-\frac{1}{2}\right) \\
&  \times B\left(  \frac{\alpha_{2}+\alpha_{3}+1}{2}-j,j+\kappa\right)
B\left(  j+\frac{1}{2},\kappa-\frac{1}{2}\right) \\
&  \times B\left(  \frac{\alpha_{3}+a_{4}+1}{2}-i-j,i+j+\kappa\right)  .
\end{align*}
We set $\alpha=0$ to find the normalizing constant, indeed%
\[
c_{\kappa}^{-1}=B\left(  2\kappa,2\kappa\right)  B\left(  \frac{1}{2}%
,\kappa\right)  ^{3}B\left(  \frac{1}{2},\kappa-\frac{1}{2}\right)  ^{2}.
\]
The stated value follows from the duplication formula. The following ratios
are typical in the calculation:%
\begin{align*}
\frac{B\left(  i+\frac{1}{2},\kappa-\frac{1}{2}\right)  }{B\left(  \frac{1}%
{2},\kappa-\frac{1}{2}\right)  }  &  =\frac{\Gamma\left(  i+\frac{1}%
{2}\right)  \Gamma\left(  \kappa-\frac{1}{2}\right)  \Gamma\left(
\kappa\right)  }{\Gamma\left(  \kappa+i\right)  \Gamma\left(  \frac{1}%
{2}\right)  \Gamma\left(  \kappa-\frac{1}{2}\right)  }=\frac{\left(  \frac
{1}{2}\right)  _{i}}{\left(  \kappa\right)  _{i}},\\
\frac{B\left(  b_{0}+\frac{1}{2}-i,i+\kappa\right)  }{B\left(  \frac{1}%
{2},\kappa\right)  }  &  =\frac{\Gamma\left(  b_{0}+\frac{1}{2}-i\right)
\Gamma\left(  i+\kappa\right)  \Gamma\left(  \kappa+\frac{1}{2}\right)
}{\Gamma\left(  b_{0}+\frac{1}{2}+\kappa\right)  \Gamma\left(  \frac{1}%
{2}\right)  \Gamma\left(  \kappa\right)  }=\frac{\left(  \kappa\right)
_{i}\left(  \frac{1}{2}\right)  _{b_{0}-i}}{\left(  \kappa+\frac{1}{2}\right)
_{b_{0}}}.
\end{align*}
Thus%
\begin{align*}
s\left(  \alpha_{1},\alpha_{2},\alpha_{3},\alpha_{4}\right)   &  =\sum
_{i=0}^{\left\lfloor \alpha_{3}/2\right\rfloor }\sum_{j=0}^{\left\lfloor
\alpha_{4}/2\right\rfloor }\binom{\alpha_{3}}{2i}\binom{\alpha_{4}}{2j}%
\frac{\left(  2\kappa\right)  _{\alpha_{1}+\alpha_{4}}\left(  2\kappa\right)
_{\alpha_{2}+\alpha_{3}}}{\left(  4\kappa\right)  _{\alpha_{1}+\alpha
_{2}+\alpha_{3}+\alpha_{4}}}\\
&  \times\frac{\left(  \frac{1}{2}\right)  _{i}\left(  \frac{1}{2}\right)
_{b_{0}-i}\left(  \frac{1}{2}\right)  _{j}\left(  \frac{1}{2}\right)
_{b_{1}-j}\left(  \frac{1}{2}\right)  _{b_{3}-i-j}\left(  \kappa\right)
_{i+j}}{\left(  \kappa+\frac{1}{2}\right)  _{b_{0}}\left(  \kappa+\frac{1}%
{2}\right)  _{b_{1}}\left(  \kappa+\frac{1}{2}\right)  _{b_{3}}}.
\end{align*}
To finish the proof, write $\binom{\alpha_{4}}{2i}=\left(  -\alpha_{4}\right)
_{2i}/\left(  2^{2i}i!\left(  \frac{1}{2}\right)  _{i}\right)  $, $\left(
\frac{1}{2}\right)  _{b_{0}-i}=\left(  -1\right)  ^{i}\left(  \frac{1}%
{2}\right)  _{b_{0}}/\left(  \frac{1}{2}-b_{0}\right)  _{i}$ and similarly for
the other terms.
\end{proof}

It is clear that the symmetry $s\left(  \alpha_{1},\alpha_{2},\alpha
_{3},\alpha_{4}\right)  =s\left(  \alpha_{2},\alpha_{1},\alpha_{4},\alpha
_{3}\right)  $ holds (as well as $s\left(  \alpha_{1},\alpha_{2},\alpha
_{3},\alpha_{4}\right)  =s\left(  \alpha_{1},\alpha_{3},\alpha_{2},\alpha
_{4}\right)  $), as implied by Lemma \ref{invmu}. Some other symmetries will
be shown later. Singular values are numerical values of $\kappa$ for which the
intertwining operator does not exist. By the general theory of singular values
\cite{DJO} specialized to $B_{2}$ they consist of $-\frac{1}{2}-\mathbb{N}%
_{0}\cup\left(  -\frac{1}{4}-\mathbb{N}_{0}\cup-\frac{3}{4}-\mathbb{N}%
_{0}\right)  $. The denominators $\left(  \kappa+\frac{1}{2}\right)  _{b_{i}}$
in $s\left(  \alpha_{1},\alpha_{2},\alpha_{3},\alpha_{4}\right)  $ correspond
to the first subset. For the second subset consider
\begin{align*}
\frac{\left(  2\kappa\right)  _{2b_{1}}\left(  2\kappa\right)  _{2b_{0}}%
}{\left(  4\kappa\right)  _{2b_{1}+2b_{0}}}  &  =2^{-2b_{1}-2b_{0}}%
\frac{\left(  2\kappa\right)  _{2b_{1}}\left(  2\kappa\right)  _{2b_{0}}%
}{\left(  2\kappa\right)  _{b_{1}+b_{0}}\left(  2\kappa+\frac{1}{2}\right)
_{b_{1}+b_{0}}}\\
&  =2^{-2b_{1}-2b_{0}}\frac{\left(  2\kappa+b_{1}+b_{0}\right)  _{b_{1}-b_{0}%
}\left(  2\kappa\right)  _{2b_{0}}}{\left(  2\kappa+\frac{1}{2}\right)
_{b_{1}+b_{0}}}%
\end{align*}
if $b_{1}\geq b_{0}$ or a similar expression if $b_{0}\geq b_{1}$.

To see some of the complexity of this integral: consider the reduced form%
\begin{gather*}
\int_{\Omega}f\left(  q_{1},q_{2}\right)  d\mu\left(  q\right)  =\\
c\int_{0}^{1}\int_{0}^{\pi}\int_{0}^{\pi}f\left(  u\cos\psi_{1},\left(
1-u\right)  \cos\psi_{2}\right)  \left(  u\left(  1-u\right)  \sin\psi_{1}%
\sin\psi_{2}\right)  ^{2\kappa-1}dud\psi_{1}d\psi_{2},
\end{gather*}
with a constant $c$ depending on $\kappa$. This is not reducible to an
ordinary double integral without the use of a nonelementary integral. Indeed,%
\begin{align*}
\int_{\Omega}f\left(  q_{1},q_{2}\right)  d\mu\left(  q\right)   &  =\int
\int_{\left\vert q_{1}\right\vert +\left\vert q_{2}\right\vert \leq1}f\left(
q_{1},q_{2}\right)  E_{\kappa}\left(  q_{1},q_{2}\right)  dq_{1}dq_{2},\\
E_{\kappa}\left(  q_{1},q_{2}\right)   &  =c\int_{\left\vert q_{1}\right\vert
}^{1-\left\vert q_{2}\right\vert }\left\{  \left(  u^{2}-q_{1}^{2}\right)
\left(  \left(  1-u\right)  ^{2}-q_{2}^{2}\right)  \right\}  ^{\kappa-1}du.
\end{align*}

\subsection{\label{sctsingl}Single sum formula and hypergeometric functions}

Even though our evaluation of $s\left(  \alpha_{1},\alpha_{2},\alpha
_{3},\alpha_{4}\right)  $ required a six-variable integral the value can be
expressed as a single sum in terms of a terminating balanced $_{4}F_{3}$
hypergeometric series.

Write
\[
s\left(  \alpha_{1},\alpha_{2},\alpha_{3},\alpha_{4}\right)  =\frac{\left(
2\kappa\right)  _{2b_{1}}\left(  2\kappa\right)  _{2b_{0}}}{\left(
4\kappa\right)  _{2b_{1}+2b_{0}}}s^{\prime}\left(  \alpha_{1},\alpha
_{2},\alpha_{3},\alpha_{4}\right)  ,
\]
then $s^{\prime}$ satisfies the recurrence%
\begin{align}
&  \alpha_{1}\alpha_{4}\left(  \kappa+\frac{1}{2}\left(  \alpha_{2}+\alpha
_{3}+1\right)  \right)  s^{\prime}\left(  \alpha_{1}-1,\alpha_{2}+1,\alpha
_{3}+1,\alpha_{4}-1\right) \label{recurs}\\
&  +\frac{1}{2}\left(  \alpha_{2}\alpha_{3}\left(  \alpha_{1}+\alpha
_{4}+1\right)  -\alpha_{1}\alpha_{4}\left(  \alpha_{2}+\alpha_{3}+1\right)
\right)  s^{\prime}\left(  \alpha_{1},\alpha_{2},\alpha_{3},\alpha_{4}\right)
\nonumber\\
&  =\alpha_{2}\alpha_{3}\left(  \kappa+\frac{1}{2}\left(  \alpha_{1}%
+\alpha_{4}+1\right)  \right)  s^{\prime}\left(  \alpha_{1}+1,\alpha
_{2}-1,\alpha_{3}-1,\alpha_{4}+1\right)  .\nonumber
\end{align}
There is a good reason why this formula appears here. The intertwining
operator will be described as a linear functional on the space of polynomials
in $q$ applied to polynomials in $x\tau\left(  q\right)  =\left(  x_{1}%
q_{1}+x_{2}q_{2},x_{1}q_{3}+x_{2}q_{4}\right)  $. The coefficient of
$x_{1}^{a+b-c}x_{2}^{c}$ in the expansion of $\left(  x_{1}q_{1}+x_{2}%
q_{2}\right)  ^{a}\left(  x_{1}q_{3}+x_{2}q_{4}\right)  ^{b}$ is the
following:%
\begin{equation}
P_{a,b}^{c}\left(  q\right)  =\sum_{i=\max\left(  0,c-a\right)  }^{\min\left(
b,c\right)  }\binom{a}{c-i}\binom{b}{i}q_{1}^{a-c+i}q_{2}^{c-i}q_{3}%
^{b-i}q_{4}^{i}. \label{defPq}%
\end{equation}
The need to integrate these polynomials motivated the examination of and
experimentation with $s\left(  a-c+i,c-i,b-i,i\right)  $ as a function of $i$.
This led to the discovery of the recurrence which in turn suggested that there
might be a single-sum form of $s$. It turns out that the proof of the
recurrence actually uses the single sum, a terminating hypergeometric series
of $_{4}F_{3}$-type (the argument looks circular, but the proof of the single
sum does not use the recurrence).

\begin{proposition}
\label{smzero}Suppose $\alpha_{4}=0$ and $\alpha_{i}=2\beta_{i}$ for $1\leq
i\leq3$, then
\begin{align*}
s^{\prime}\left(  2\beta_{1},2\beta_{2},2\beta_{3},0\right)   &
=\frac{\left(  \frac{1}{2}\right)  _{\beta_{1}}\left(  \frac{1}{2}\right)
_{\beta_{2}}\left(  \frac{1}{2}\right)  _{\beta_{3}}}{\left(  \kappa+\frac
{1}{2}\right)  _{\beta_{1}}\left(  \kappa+\frac{1}{2}\right)  _{\beta_{2}%
}\left(  \kappa+\frac{1}{2}\right)  _{\beta_{3}}}\\
&  =s^{\prime}\left(  0,2\beta_{2},2\beta_{3},2\beta_{1}\right)  ,
\end{align*}
and the value is symmetric in $\beta_{1},\beta_{2},\beta_{3}$.
\end{proposition}

\begin{proof}
Indeed the sum in $s^{\prime}\left(  2\beta_{1},2\beta_{2},2\beta
_{3},0\right)  $ equals
\begin{align*}
&  \sum_{j=0}^{\beta_{3}}\frac{\left(  -\beta_{3}\right)  _{j}\left(  \frac
{1}{2}-\beta_{3}\right)  _{j}\left(  \kappa\right)  _{j}}{j!\left(  \frac
{1}{2}-\beta_{2}-\beta_{3}\right)  _{j}\left(  \frac{1}{2}-\beta_{3}\right)
_{j}}\\
&  =\frac{\left(  \frac{1}{2}-\beta_{2}-\beta_{3}-\kappa\right)  _{\beta_{3}}%
}{\left(  \frac{1}{2}-\beta_{2}-\beta_{3}\right)  _{\beta_{3}}}=\frac{\left(
\kappa+\beta_{2}+\frac{1}{2}\right)  _{\beta_{3}}}{\left(  \frac{1}{2}%
+\beta_{2}\right)  _{\beta_{3}}}\\
&  =\frac{\left(  \kappa+\frac{1}{2}\right)  _{\beta_{2}+\beta_{3}}\left(
\frac{1}{2}\right)  _{\beta_{2}}}{\left(  \kappa+\frac{1}{2}\right)
_{\beta_{2}}\left(  \frac{1}{2}\right)  _{\beta_{2}+\beta_{3}}},
\end{align*}
(by the Chu-Vandermonde sum) and this proves the first equation. The second
equation follows from the symmetry of Lemma \ref{invmu}. The sum in
$s^{\prime}\left(  0,2\beta_{2},2\beta_{3},2\beta_{1}\right)  $ can also be
found directly, first summing over $0\leq i\leq\beta_{1}$ and using similar
arguments as in the first equation.
\end{proof}

\begin{theorem}
Suppose $\alpha\in\mathbb{N}_{0}^{4}$ and $\alpha_{1}\equiv\alpha_{2}%
\equiv\alpha_{3}\equiv\alpha_{4}\operatorname{mod}2$. Let $b_{i}=\left(
\alpha_{i}+\alpha_{4}\right)  /2,1\leq i\leq3$, then%
\begin{align}
s^{\prime}\left(  \alpha_{1},\alpha_{2},\alpha_{3},\alpha_{4}\right)   &
=\frac{\left(  \frac{1}{2}\right)  _{b_{1}}\left(  \frac{1}{2}\right)
_{b_{2}}\left(  \frac{1}{2}\right)  _{b_{3}}}{\left(  \kappa+\frac{1}%
{2}\right)  _{b_{1}}\left(  \kappa+\frac{1}{2}\right)  _{b_{2}}\left(
\kappa+\frac{1}{2}\right)  _{b_{3}}}\times\label{sngl}\\
&  \sum_{i=0}^{\left\lfloor \alpha_{4}/2\right\rfloor }\frac{\left(
-\frac{\alpha_{4}}{2}\right)  _{i}\left(  \frac{1-\alpha_{4}}{2}\right)
_{i}\left(  \kappa\right)  _{i}\left(  -\kappa-b_{1}-b_{0}\right)  _{i}%
}{i!\left(  \frac{1}{2}-b_{1}\right)  _{i}\left(  \frac{1}{2}-b_{2}\right)
_{i}\left(  \frac{1}{2}-b_{3}\right)  _{i}}.\nonumber
\end{align}
This is a terminating balanced $_{4}F_{3}$ series.
\end{theorem}

The proof of this key result is in Section 3. To say that the series is
balanced means that the sum of the denominator parameters equals one plus the
sum of the numerator parameters; the property is also called
Saalsch\"{u}tzian. The particular choice of parameters for this $_{4}F_{3}$
series will appear often in the sequel and so we make the following:

\begin{definition}
\label{defF}For $n\in\mathbb{N}_{0}$ and free parameters $u,v_{1},v_{2},v_{3}$
let%
\[
F\left(  n;u,v_{1},v_{2},v_{3}\right)  =~_{4}F_{3}\left(
\genfrac{}{}{0pt}{}{-\frac{n}{2},\frac{1-n}{2},u,-u-v_{1}-v_{2}-v_{3}%
}{\frac{1}{2}-v_{1}-n,\frac{1}{2}-v_{2},\frac{1}{2}-v_{3}}%
;1\right)  .
\]
This is a balanced terminating hypergeometric function.
\end{definition}

This is not the generic balanced $_{4}F_{3}$-series (which has 6 free
parameters) because of the parameters $\left(  -\frac{n}{2},\frac{1-n}%
{2}\right)  $. With this notation the Theorem can be restated as:

\begin{theorem}
For $\alpha_{1}\equiv\alpha_{2}\equiv\alpha_{3}\equiv\alpha_{4}%
\operatorname{mod}2$ the following single-sum expression is valid:%
\begin{align*}
s\left(  \alpha\right)   &  =\frac{\left(  2\kappa\right)  _{\alpha_{1}%
+\alpha_{4}}\left(  2\kappa\right)  _{\alpha_{2}+\alpha_{3}}}{\left(
4\kappa\right)  _{\left\vert \alpha\right\vert }}\frac{\left(  \frac{1}%
{2}\right)  _{b_{1}}\left(  \frac{1}{2}\right)  _{b_{2}}\left(  \frac{1}%
{2}\right)  _{b_{3}}}{\left(  \kappa+\frac{1}{2}\right)  _{b_{1}}\left(
\kappa+\frac{1}{2}\right)  _{b_{2}}\left(  \kappa+\frac{1}{2}\right)  _{b_{3}%
}}\\
&  \times F\left(  \alpha_{4};\kappa,\frac{1}{2}\left(  \alpha_{1}-\alpha
_{4}\right)  ,\frac{1}{2}\left(  \alpha_{2}+\alpha_{4}\right)  ,\frac{1}%
{2}\left(  \alpha_{3}+\alpha_{4}\right)  \right)  .
\end{align*}

\end{theorem}

\subsection{Contents}

In Section 2 we describe the representation-theoretic implications of the
invariance conditions, present an overview of the method of proving that a
given linear functional on polynomials in $q$ produces the intertwining
operator and then give the actual proof. Some ingredients of the proofs depend
on contiguity relations for the function $F$. Also there is a purely integral
form of $V$ not involving the degree of polynomials.

Section 3 contains the proof of the single-sum result, which relies on
classical transformations of hypergeometric series, and the
(computer-assisted) proofs of the required contiguity relations of $F$. There
are closing comments in Section 4.

\section{The Intertwining Operator}

\subsection{\label{sctB2}Invariants and representations of $B_{2}$}

The irreducible representations of $B_{2}$ are realized in the space of
polynomials $\mathbb{R}\left[  x_{1},x_{2}\right]  $ as follows: there is one
of degree 2, which is realized in each space of polynomials homogeneous of
degree $2n-1$ as $n$ isomorphic copies%
\[
\left\{  c_{0}x_{1}^{2n-1-j}x_{2}^{j}+c_{1}x_{1}^{j}x_{2}^{2n-1-j}:0\leq j\leq
n-1\right\}  ,
\]
for $n=1,2,\ldots$; there are four one-dimensional representations:

\begin{enumerate}
\item invariants: $x_{1}^{2n-2j}x_{2}^{2j}+x_{1}^{2j}x_{2}^{2n-2j},0\leq j\leq
n$;

\item determinant: $x_{1}^{2n-1-2j}x_{2}^{2j+1}-x_{1}^{2j+1}x_{2}%
^{2n-1-2j},0\leq j\leq\left\lfloor \frac{n}{2}\right\rfloor -1$ and $2n\geq4$;

\item type 1: $x_{1}^{2n-1-2j}x_{2}^{2j+1}+x_{1}^{2j+1}x_{2}^{2n-1-2j},0\leq
j\leq\left\lfloor \frac{n-1}{2}\right\rfloor $ and $2n\geq2$;

\item type 2: $x_{1}^{2n-2j}x_{2}^{2j}-x_{1}^{2j}x_{2}^{2n-2j},0\leq
j\leq\left\lfloor \frac{n-1}{2}\right\rfloor $ and $2n\geq2$.
\end{enumerate}

The reason for naming types is that $\sigma_{i}$ acts as multiplication by
$-1$ on type $\#i$ for $i=1,2$. Consider polynomials in $q$ invariant under
$q\mapsto q\lambda\left(  w\right)  \rho\left(  w\right)  ,w\in B_{2}$. It is
easy to compute the Poincar\'{e} series for the ring of invariants (graded by
degree), namely $\frac{1+z^{3}}{\left(  1-z\right)  \left(  1-z^{2}\right)
^{3}}$ and then describe the ring as
\[
\left(  1\oplus\left(  q_{1}-q_{4}\right)  \left(  q_{2}^{2}-q_{3}^{2}\right)
\right)  \mathbb{R}\left[  q_{1}+q_{4},q_{1}^{2}+q_{4}^{2},q_{2}^{2}+q_{3}%
^{2},q_{2}q_{3}\right]  .
\]
The invariant of degree 1 is associated with the two-dimensional
representation of $B_{2}$ and we set
\[
g_{0}\left(  q\right)  =2\left(  q_{1}+q_{4}\right)  .
\]
There are 4 linearly independent invariants of degree 2, the two-sided
$\sum_{i=1}^{4}q_{i}^{2}$, and%
\begin{align*}
g_{1}\left(  q\right)   &  =q_{1}q_{4}+q_{2}q_{3},\\
g_{2}\left(  q\right)   &  =\frac{1}{2}\left(  q_{1}^{2}-q_{2}^{2}-q_{3}%
^{2}+q_{4}^{2}\right)  ,\\
g_{3}\left(  q\right)   &  =q_{1}q_{4}-q_{2}q_{3}.
\end{align*}
Then%
\begin{align*}
g_{1}\left(  q\lambda\left(  \sigma_{1}\right)  \right)   &  =-g_{1}\left(
q\right)  ,g_{2}\left(  q\lambda\left(  \sigma_{1}\right)  \right)
=g_{2}\left(  q\right)  ,g_{3}\left(  q\lambda\left(  \sigma_{1}\right)
\right)  =-g_{3}\left(  q\right)  ,\\
g_{1}\left(  q\lambda\left(  \sigma_{2}\right)  \right)   &  =g_{1}\left(
q\right)  ,g_{2}\left(  q\lambda\left(  \sigma_{2}\right)  \right)
=-g_{2}\left(  q\right)  ,g_{3}\left(  q\lambda\left(  \sigma_{2}\right)
\right)  =-g_{3}\left(  q\right)  .
\end{align*}
Further $g_{i}\left(  1,0,0,1\right)  =1$ for $1\leq i\leq3$ (a plausible
normalization). By simple orthogonality arguments we see that if $f\left(
x\right)  $ is a polynomial of even degree then $\int_{\Omega}f\left(
x\tau\left(  q\right)  \right)  g_{i}\left(  q\right)  d\mu\left(  q\right)
\neq0$ only if $f\left(  x\right)  $ has a nonzero component of determinant
type for $i=3$, or type $i$ for $i=1,2$. Experimentation quickly showed that
the formula $Vf\left(  x\right)  =\int_{\Omega}f\left(  x\tau\left(  q\right)
\right)  g_{3}\left(  q\right)  d\mu\left(  q\right)  $ appeared to be valid
for the determinant type, but the similar attempt failed for types 1 and 2 (in
fact, no polynomial in $q$ of degree less than twelve with the correct
behavior under the $B_{2}$-action works). We comment on this quandary in
Section 4. We set
\[
\partial_{i}=\frac{\partial}{\partial q_{i}},1\leq i\leq4.
\]
The formula for $V$ comes from applying the adjoint action $L:=\mathrm{ad}%
\left(  \frac{1}{2}\sum_{i=1}^{4}\partial_{i}^{2}\right)  $ to (multiplication
by) the $g_{i}$ (recall $\left(  \mathrm{ad}\left(  A\right)  B\right)
f:=\left(  AB-BA\right)  f$ for operators $A,B$ on polynomials $f$). Indeed%
\begin{align*}
L\left(  g_{1}\right)   &  =q_{1}\partial_{4}+q_{4}\partial_{1}+q_{2}%
\partial_{3}+q_{3}\partial_{2},\\
L\left(  g_{2}\right)   &  =q_{1}\partial_{1}-q_{2}\partial_{2}-q_{3}%
\partial_{3}+q_{4}\partial_{4},\\
L\left(  g_{3}\right)   &  =q_{1}\partial_{4}+q_{4}\partial_{1}-q_{2}%
\partial_{3}-q_{3}\partial_{2}.
\end{align*}
These operators have the same invariance properties as the respective
polynomials $g_{i}$. Set $D_{0}:=L\left(  g_{1}+g_{2}\right)  =\left(
q_{1}+q_{4}\right)  \left(  \partial_{1}+\partial_{4}\right)  -\left(
q_{2}-q_{3}\right)  \left(  \partial_{2}-\partial_{3}\right)  $ and
$D_{3}:=L\left(  g_{3}\right)  $.

The formula stated in the introduction can be given as a pure integral with no
derivatives for $\kappa>\frac{3}{2}$. However another variable of integration
occurs (a total of seven!). Let $\omega$ denote the integration operator:%
\[
\omega f\left(  x\right)  :=\int_{0}^{1}f\left(  tx\right)  t^{4\kappa-1}dt,
\]
thus $\omega f\left(  x\right)  =\frac{1}{4\kappa+n}f\left(  x\right)  $ when
$f$ is homogeneous of degree $n$.

\begin{theorem}
\label{Vint}Let $\kappa>\frac{3}{2}$ and suppose $f$ is sufficiently smooth on
a ball $B_{R}=\left\{  x\in\mathbb{R}^{2}:\left\Vert x\right\Vert
_{2}<R\right\}  $ for some $R>0$, then for $x\in B_{R}$
\begin{align*}
Vf\left(  x\right)   &  =\int_{\Omega}f\left(  x\tau\left(  q\right)  \right)
\left(  1+g_{0}\left(  q\right)  +g_{3}\left(  q\right)  \right)  d\mu\left(
q\right) \\
&  +\left(  2\kappa-3\right)  \int_{\Omega}\omega f\left(  x\tau\left(
q\right)  \right)  \widetilde{g}_{0}\left(  q,u,\theta,\phi_{1},\phi_{2}%
,\psi_{1},\psi_{2}\right)  d\mu\left(  q\right)  ,\\
\widetilde{g}_{0}  &  =\frac{1}{1+\cos\theta}\left(  \frac{q_{1}+q_{4}}%
{u\sin\phi_{1}\sin\psi_{1}}\right)  ^{2}-\frac{1}{1-\cos\theta}\left(
\frac{q_{2}-q_{3}}{\left(  1-u\right)  \sin\phi_{2}\sin\psi_{2}}\right)  ^{2}.
\end{align*}

\end{theorem}

If $f\left(  x\right)  $ is homogeneous of degree $n$ then $\int_{\Omega
}f\left(  x\tau\left(  q\right)  \right)  g_{3}\left(  q\right)  d\mu\left(
q\right)  $ is zero when $n$ is odd and equals $\frac{1}{8\kappa+n}%
\int_{\Omega}D_{3}f\left(  x\tau\left(  q\right)  \right)  d\mu\left(
q\right)  $ when $n$ is even (the proof is in the next subsection). Also we
will prove the two equations (\ref{oddfV}) and (\ref{evenfV}). The factor
$\frac{1}{4\kappa+n}$ in (\ref{evenfV}) can be replaced by the integral
operator $\omega$. The following uses an integration by parts to replace the
differential operator $D_{0}$ by an integral, the last ingredient of the
formula in the theorem.

\begin{lemma}
For $\kappa>\frac{3}{2}$ and a smooth function $h\left(  q\right)  $%
\begin{gather*}
\int_{\Omega}D_{0}h\left(  q\right)  d\mu\left(  q\right)  =\left(
2\kappa-3\right)  \times\\
\int_{\Omega}h\left(  q\right)  \left\{  \frac{1}{1+\cos\theta}\left(
\frac{q_{1}+q_{4}}{u\sin\phi_{1}\sin\psi_{1}}\right)  ^{2}-\frac{1}%
{1-\cos\theta}\left(  \frac{q_{2}-q_{3}}{\left(  1-u\right)  \sin\phi_{2}%
\sin\psi_{2}}\right)  ^{2}\right\}  d\mu\left(  q\right)
\end{gather*}

\end{lemma}

\begin{proof}
Use the notation and change-of-variable from Lemma \ref{invmu} to set up the
integration by parts. The measure is $d\mu=c_{\kappa}\left(  z_{1}%
z_{2}\right)  ^{\kappa-3/2}\left(  \sin\theta\right)  ^{2\kappa-3}dq_{1}%
dq_{2}dq_{3}dq_{4}dud\theta$ and the region of integration $\Omega_{1}$ is
implicitly defined by $z_{1}\geq0,z_{2}\geq0,0\leq\theta\leq\pi,0\leq u\leq1$.
If $w\left(  q\right)  $ vanishes on the boundary of the bounded domain $Q$ in
$\mathbb{R}^{4}$ then $\int_{Q}D_{0}h\left(  q\right)  w\left(  q\right)
dq=-\int_{Q}h\left(  q\right)  D_{0}w\left(  q\right)  dq,$ (although
$\int_{Q}q_{j}\frac{\partial}{\partial q_{j}}h\left(  q\right)  w\left(
q\right)  dq=-\int_{Q}h\left(  q\right)  \left(  1+q_{j}\frac{\partial
}{\partial q_{j}}\right)  w\left(  q\right)  dq$ for $1\leq j\leq4$ the terms
$\int_{Q}h\left(  q\right)  w\left(  q\right)  dq$ cancel out). We find
\begin{align*}
D_{0}z_{1}  &  =-2\frac{\left(  1-\cos\theta\right)  \left(  q_{1}%
+q_{4}\right)  ^{2}}{\sin^{2}\theta}=-2\frac{\left(  q_{1}+q_{4}\right)  ^{2}%
}{1+\cos\theta},\\
D_{0}z_{2}  &  =2\frac{\left(  1+\cos\theta\right)  \left(  q_{2}%
-q_{3}\right)  ^{2}}{\sin^{2}\theta}=2\frac{\left(  q_{2}-q_{3}\right)  ^{2}%
}{1-\cos\theta},
\end{align*}
and $D_{0}\left(  z_{1}z_{2}\right)  ^{\kappa-3/2}=\left(  \kappa-\frac{3}%
{2}\right)  \left(  \frac{D_{0}z_{1}}{z_{1}}+\frac{D_{0}z_{2}}{z_{2}}\right)
\left(  z_{1}z_{2}\right)  ^{\kappa-3/2}$. Then change back to the original
variables $\phi_{1},\phi_{2},\psi_{1},\psi_{2},u,\theta$ to complete the proof.
\end{proof}

\subsection{Proof of the intertwining property}

Suppose $\xi$ is a linear functional on polynomials in $q=\left(  q_{1}%
,q_{2},q_{3},q_{4}\right)  $, and define an operator on polynomials in
$x=\left(  x_{1},x_{2}\right)  $ by $V_{1}f\left(  x\right)  =\xi f\left(
x_{1}q_{1}+x_{2}q_{2},x_{1}q_{3}+x_{2}q_{4}\right)  $. What needs to be done
to show $V_{1}=V$, the intertwining operator? The group invariance requires
$\xi p\left(  q\right)  =\xi p\left(  q_{4},q_{3},q_{2},q_{1}\right)  =\xi
p\left(  q_{1},-q_{2},-q_{3},q_{4}\right)  $.

\begin{definition}
Let $\xi_{0}$ denote the functional defined by $q^{\alpha}\mapsto s\left(
\alpha\right)  $. For a homogeneous polynomial $p\left(  q\right)  $ in $q$ of
degree $n$ define the functional $\xi$ by:%
\begin{align*}
\xi\left(  p\right)   &  =\xi_{0}\left(  g_{0}p\right)  ,\text{ for }n\text{
odd,}\\
\xi\left(  p\right)   &  =\xi_{0}\left(  \left(  1+\frac{1}{4\kappa+n}%
D_{0}+\frac{1}{8\kappa+n}D_{3}\right)  p\left(  q\right)  \right)  ,\text{ for
}n\text{ even.}%
\end{align*}

\end{definition}

Using \cite[Prop. 1.3]{D3} we get the homogeneous component of the criterion:%
\begin{align}
&  \left(  n+1\right)  \left(  \left\langle x,y\right\rangle \xi\left(
\left\langle x\tau\left(  q\right)  ,y\right\rangle ^{n}\right)  -\xi\left(
\left\langle x\tau\left(  q\right)  ,y\right\rangle ^{n+1}\right)  \right)
\label{condV}\\
&  =\kappa\sum_{i=1}^{4}\left(  \xi\left(  \left\langle x\tau\left(  q\right)
,y\right\rangle ^{n+1}\right)  -\xi\left(  \left\langle x\sigma_{i}\tau\left(
q\right)  ,y\right\rangle ^{n+1}\right)  \right)  ,\nonumber
\end{align}
for $n=0,1,2,\ldots$, where $\left\{  \sigma_{i}\right\}  $ is the set of
reflections $\left\{  \sigma_{1},\sigma_{2},\sigma_{2}\sigma_{1}\sigma
_{2},\sigma_{1}\sigma_{2}\sigma_{1}\right\}  $, and $\sigma_{i}\tau\left(
q\right)  =\tau\left(  q\lambda\left(  \sigma_{i}\right)  \right)  $. Rewrite
the criterion as:%
\begin{gather*}
\xi\left(  \left(  \partial_{1}+\partial_{4}\right)  \left\langle x\tau\left(
q\right)  ,y\right\rangle ^{n+1}\right)  -\left(  n+1+4\kappa\right)
\xi\left(  \left\langle x\tau\left(  q\right)  ,y\right\rangle ^{n+1}\right)
\\
+\kappa\sum_{i=1}^{4}\xi\left(  \left\langle x\tau\left(  q\lambda\left(
\sigma_{i}\right)  \right)  ,y\right\rangle ^{n+1}\right)  =0
\end{gather*}
Note $\left(  \partial_{1}+\partial_{4}\right)  \left\langle x\tau\left(
q\right)  ,y\right\rangle =\left\langle x,y\right\rangle $. We will prove
$\xi$ satisfies the criterion and also that $\xi_{0}\left(  D_{3}p\left(
q\right)  \right)  =\left(  8\kappa+n\right)  \xi_{0}\left(  g_{3}p\right)  $
when $p$ is of degree $n$ (both sides vanish when $n$ is odd). We have%
\begin{align*}
q\lambda\left(  \sigma_{1}\right)   &  =\left(  -q_{1},q_{2},-q_{3}%
,q_{4}\right)  ,\\
q\lambda\left(  \sigma_{2}\right)   &  =\left(  q_{2},q_{1},q_{4}%
,q_{3}\right)  ,\\
q\lambda\left(  \sigma_{1}\sigma_{2}\sigma_{1}\right)   &  =\left(
-q_{2},-q_{1},-q_{4},-q_{3}\right)  ,\\
q\lambda\left(  \sigma_{2}\sigma_{1}\sigma_{2}\right)   &  =\left(
q_{1},-q_{2},q_{3},-q_{4}\right)  .
\end{align*}
Then
\begin{align*}
\sum_{i=1}^{4}g_{j}\left(  q\lambda\left(  \sigma_{i}\right)  \right)   &
=0,j=0,1,2,\\
\sum_{i=1}^{4}g_{3}\left(  q\lambda\left(  \sigma_{i}\right)  \right)   &
=-4g_{3}\left(  q\right)  .
\end{align*}
The corresponding differential operators $L\left(  g_{j}\right)  $ satisfy
similar equations for $j=1,2,3$ because $\sum_{i=1}^{4}\partial_{i}^{2}$ is
two-sided invariant, that is, $\sum_{i=1}^{4}L\left(  g_{j}\right)  p\left(
q\lambda\left(  \sigma_{i}\right)  \right)  =0$ for $j=1,2$ and $\sum
_{i=1}^{4}L\left(  g_{3}\right)  p\left(  q\lambda\left(  \sigma_{i}\right)
\right)  =-4L\left(  g_{3}\right)  p\left(  q\right)  $.

The proof of the criterion is easy when $n$ is odd.

\begin{proposition}
Suppose $n$ is odd and $p\left(  q\right)  $ is homogeneous of degree $n+1$,
then
\begin{gather*}
\xi_{0}\left(  g_{0}\left(  q\right)  \left(  \partial_{1}+\partial
_{4}\right)  p\left(  q\right)  \right) \\
-\left(  n+1+4\kappa\right)  \xi_{0}\left(  \left(  1+\frac{1}{4\kappa
+n+1}D_{0}+\frac{1}{8\kappa+n+1}D_{3}\right)  p\left(  q\right)  \right) \\
+\kappa\xi_{0}\left(  \left(  4-\frac{4}{8\kappa+n+1}D_{3}\right)  p\left(
q\right)  \right)  =0.
\end{gather*}

\end{proposition}

\begin{proof}
In fact, the left side simplifies to
\[
\xi_{0}\left(  \left(  2\left(  q_{1}+q_{4}\right)  \left(  \partial
_{1}+\partial_{4}\right)  -\left(  n+1\right)  -D_{0}-D_{3}\right)  p\left(
q\right)  \right)  .
\]
Replace $n+1$ by the Euler operator $\sum_{i=1}^{4}q_{i}\partial_{i}$ then the
expression becomes identically zero.
\end{proof}

\begin{corollary}
The criterion (\ref{condV}) is satisfied for odd $n.$
\end{corollary}

\begin{proof}
Set $p\left(  q\right)  =\left\langle x\tau\left(  q\right)  ,y\right\rangle
^{n+1}$ in the Theorem.
\end{proof}

For the odd degree case ($n$ even) replace $n$ by $2n$, then the criterion
becomes
\begin{equation}
\left(  4\kappa+2n+1\right)  \xi_{0}\left(  \left\langle x\tau\left(
q\right)  ,y\right\rangle ^{2n+1}g_{0}\left(  q\right)  \right)  -\xi\left(
\left(  \partial_{1}+\partial_{4}\right)  \left\langle x\tau\left(  q\right)
,y\right\rangle ^{2n+1}\right)  =0. \label{oddeqn}%
\end{equation}
This equation can be restated as
\[
\left(  4\kappa+2n+1\right)  K_{2n+1}\left(  x,y\right)  =\left\langle
x,y\right\rangle K_{2n}\left(  x,y\right)  .
\]
Recall the definition of $P_{a,b}^{c}\left(  q\right)  $ from (\ref{defPq}).
Thus%
\[
\left\langle x\tau\left(  q\right)  ,y\right\rangle ^{n}=\sum_{i=0}^{n}%
\binom{n}{i}y_{1}^{n-i}y_{2}^{i}\sum_{c=0}^{n}P_{n-i,i}^{c}\left(  q\right)
x_{1}^{n-c}x_{2}^{c},
\]
for $n\in\mathbb{N}_{0}$. So it suffices to prove identities involving $\xi$
and $\left\langle x\tau\left(  q\right)  ,y\right\rangle ^{n}$ for the
polynomials $P_{a,b}^{c}$. There are two immediate consequences of the
invariance properties of $\xi$.%
\begin{align*}
P_{a,b}^{c}\left(  q\lambda\left(  \sigma_{1}\right)  \rho\left(  \sigma
_{1}\right)  \right)   &  =\left(  -1\right)  ^{b+c}P_{a,b}^{c}\left(
q\right)  ,\\
\xi\left(  P_{a,b}^{c}\right)   &  \neq0\Longrightarrow b\equiv
c\operatorname{mod}2,
\end{align*}
and%
\begin{align*}
P_{a,b}^{c}\left(  q\lambda\left(  \sigma_{2}\right)  \rho\left(  \sigma
_{2}\right)  \right)   &  =P_{b,a}^{a+b-c}\left(  q\right)  ,\\
\xi\left(  P_{a,b}^{c}\right)   &  =\xi\left(  P_{b,a}^{a+b-c}\right)  .
\end{align*}
With the aim of applying condition (\ref{oddeqn}) to $P_{a,b}^{c}$ with
$a+b=2n+1$ we can assume that $a$ is odd, $b$ and $c$ are even. This implies
$\xi_{0}\left(  q_{4}P_{a,b}^{c}\left(  q\right)  \right)  =0$ (the typical
monomial is $q_{1}^{a-c+i}q_{2}^{c-i}q_{3}^{b-i}q_{4}^{i+1}$; the parities of
the exponents of\ $q_{3}$ and $q_{4}$ are opposite).

In the evaluation of $\xi_{0}\left(  g_{0}P_{a,b}^{c}\right)  $ and
$\xi\left(  \left(  \partial_{1}+\partial_{4}\right)  P_{a,b}^{c}\right)  $
there are several vanishing terms:%
\begin{align*}
\xi_{0}\left(  q_{4}P_{a,b}^{c}\right)   &  =\xi_{0}\left(  \partial
_{4}P_{a,b}^{c}\right)  =\xi_{0}\left(  q_{i}\partial_{i}\partial_{4}%
P_{a,b}^{c}\right)  =0,1\leq i\leq4,\\
\xi_{0}\left(  q_{4}\partial_{1}^{2}P_{a,b}^{c}\right)   &  =\xi_{0}\left(
q_{1}\partial_{4}\partial_{1}P_{a,b}^{c}\right)  =\xi_{0}\left(  q_{2}%
\partial_{3}\partial_{1}P_{a,b}^{c}\right)  =\xi_{0}\left(  q_{3}\partial
_{2}\partial_{1}P_{a,b}^{c}\right)  =0.
\end{align*}
Thus%
\begin{align*}
\xi_{0}\left(  D_{0}\left(  \partial_{1}+\partial_{4}\right)  P_{a,b}%
^{c}\right)   &  =\xi_{0}\left(  L\left(  g_{2}\right)  \partial_{1}%
P_{a,b}^{c}\right)  +\xi_{0}\left(  L\left(  g_{1}\right)  \partial_{4}%
P_{a,b}^{c}\right)  ,\\
\xi_{0}\left(  D_{3}\left(  \partial_{1}+\partial_{4}\right)  P_{a,b}%
^{c}\right)   &  =\xi_{0}\left(  D_{3}\partial_{4}P_{a,b}^{c}\right)  .
\end{align*}
The required identity is%
\begin{align}
&  \left(  4\kappa+2n+1\right)  \xi_{0}\left(  2q_{1}P_{a,b}^{c}\right)
-\xi_{0}\left(  \left(  1+\frac{1}{4\kappa+2n}L\left(  g_{2}\right)  \right)
\partial_{1}P_{a,b}^{c}\right) \label{oddP}\\
&  -\frac{1}{\left(  2\kappa+n\right)  \left(  4\kappa+n\right)  }\xi
_{0}\left(  \left\{  \left(  3\kappa+n\right)  \left(  q_{4}\partial_{1}%
+q_{1}\partial_{4}\right)  +\kappa\left(  q_{3}\partial_{2}+q_{2}\partial
_{3}\right)  \right\}  \partial_{4}P_{a,b}^{c}\right) \nonumber\\
&  =0.\nonumber
\end{align}
We will prove this by summing over the monomials in $P_{a,b}^{c}$, that is, we
replace $P_{a,b}^{c}$ in the left side by $q_{1}^{a-c+i}q_{2}^{c-i}q_{3}%
^{b-i}q_{4}^{i}$ and evaluate $\xi_{0}$ (in terms of $s\left(  \alpha\right)
$). The partial sums are found explicitly by use of contiguity relations for
the $_{4}F_{3}$-type function $F$.

Set $a=2a_{1}+1+2a_{3},b=2a_{2},c=2a_{3}$ and so $n=a_{1}+a_{2}+a_{3}$ and
assume for now that $a_{1}\geq0$, that is, $a>c$. Evaluate the left side for
the monomial $q^{\alpha}$ with $\alpha=\left(  2a_{1}+1+i,2a_{3}%
-i,2a_{2}-i,i\right)  $ (note $\left\vert \alpha\right\vert =2n+1$). One of
the terms simplifies:%
\begin{gather*}
\xi_{0}\left(  \left(  1+\frac{1}{4\kappa+2n}L\left(  g_{2}\right)  \right)
\partial_{1}q^{\alpha}\right)  =\\
\left(  2a_{1}+1+i\right)  \left(  1+\frac{2a_{1}-2a_{2}-2a_{3}-4i}%
{4\kappa+2n}\right)  s\left(  2a_{1}+i,2a_{3}-i,2a_{2}-i,i\right) \\
=2\left(  2a_{1}+1+i\right)  \frac{\left(  \kappa+a_{1}-i\right)  }{2\kappa
+n}s\left(  2a_{1}+i,2a_{3}-i,2a_{2}-i,i\right)  .
\end{gather*}

To remove some common factors we divide by $s\left(  2a_{1}+2,2a_{3}%
,2a_{2},0\right)  $ (see Proposition \ref{smzero} for the evaluation) and
denote the result by $t_{i}$. The following expression is a linear combination
of $F$ values with simple coefficients
\begin{gather*}
s\left(  2a_{1}+2,2a_{3},2a_{2},0\right)  t_{i}=2\left(  4\kappa+2n+1\right)
s\left(  2a_{1}+2+i,2a_{3}-i,2a_{2}-i,i\right) \\
-2\frac{\left(  2a_{1}+1+i\right)  \left(  \kappa+a_{1}+i\right)  }{2\kappa
+n}s\left(  2a_{1}+i,2a_{3}-i,2a_{2}-i,i\right) \\
-\frac{i\left(  3\kappa+n\right)  }{\left(  2\kappa+n\right)  \left(
4\kappa+n\right)  }\{\left(  2a_{1}+1+i\right)  s\left(  2a_{1}+i,2a_{3}%
-i,2a_{2}-i,i\right) \\
+\left(  i-1\right)  s\left(  2a_{1}+2+i,2a_{3}-i,2a_{2}-i,i-2\right)  \}\\
-\frac{i\kappa}{\left(  2\kappa+n\right)  \left(  4\kappa+n\right)  }\{\left(
2a_{3}-i\right)  s\left(  2a_{1}+1+i,2a_{3}-i-1,2a_{2}-i+1,i-1\right) \\
+\left(  2a_{2}-i\right)  s\left(  2a_{1}+1+i,2a_{3}-i+1,2a_{2}%
-i-1,i-1\right)  \}.
\end{gather*}
Then $t_{0}=0$ (directly). Rewrite $\binom{a}{c-i}\binom{b}{i}=\binom{a}%
{c}\frac{\left(  -c\right)  _{i}\left(  -b\right)  _{i}}{i!\left(
a-c+1\right)  _{i}}=\binom{a}{c}\frac{\left(  -2a_{3}\right)  _{i}\left(
-2a_{2}\right)  _{i}}{i!\left(  2a_{1}+2\right)  _{i}}$. The proof of the
following is in Section \ref{sctcont}.

\begin{theorem}
\label{big1}Suppose $a_{1}\geq0$ and $m\geq1$, then%
\begin{gather*}
\sum_{i=1}^{m}\frac{\left(  -2a_{3}\right)  _{i}\left(  -2a_{2}\right)  _{i}%
}{i!\left(  2a_{1}+2\right)  _{i}}t_{i}=\frac{2^{2m+3}\kappa a_{2}a_{3}\left(
2-2a_{2}\right)  _{m-1}\left(  2-2a_{3}\right)  _{m-1}\left(  \kappa
+a_{1}+1\right)  _{m}}{\left(  m-1\right)  !\left(  -2\kappa-2a_{2}%
-2a_{3}+1\right)  _{2m}}\\
\times\frac{\left(  a_{1}+\frac{3}{2}\right)  _{m-1}\left(  4\kappa
+2n+1\right)  }{\left(  2a_{1}+2\right)  _{m}\left(  4\kappa+n\right)
}F\left(  m-1;\kappa+1,a_{1}+1,a_{2}-1,a_{3}-1\right)  .
\end{gather*}

\end{theorem}

\begin{corollary}
The identity (\ref{oddP}) is valid and criterion (\ref{condV}) is satisfied
for even $n.$
\end{corollary}

\begin{proof}
The factor $\left(  2-2a_{2}\right)  _{m-1}\left(  2-2a_{3}\right)  _{m-1}$
vanishes for $m\geq2a_{2}$ or $m\geq2a_{3}$. Multiply both sides by
$\binom{2a_{1}+2a_{3}+1}{2a_{3}}$, the left side becomes%
\[
\sum_{i=0}^{2a_{3}}\frac{\left(  -2a_{2}\right)  _{i}\left(  -2a_{1}%
-2a_{3}-1\right)  _{2a_{3}-i}}{i!\left(  2a_{3}-i\right)  !}t_{i}.
\]
If $a_{1}\geq0$ set $m=\min\left(  2a_{2},2a_{3}\right)  $. The poles of $F$
(as rational function of $a_{2},a_{3}$) occur in a subset of $\frac{1}%
{2}+\mathbb{N}_{0}$. Next we use a weak form of analytic continuation to apply
the Theorem to the case $a_{1}<0$. The terms in the left side with
$2a_{3}-i-1+\left(  -2a_{1}-2a_{3}-1\right)  \geq0$ vanish, that is, for
$i<-2a_{1}-1$. Further $\left(  2a_{1}+1\right)  !\times\left(  2a_{1}%
+2\right)  _{m}=\left(  2a_{1}+m+1\right)  !$). Since $1/\left(
2a_{1}+m+1\right)  !$ is entire for fixed $m$ the identity (left side minus
right side, for fixed $a_{2},a_{3},m$) can be considered as a meromorphic
function of $a_{1}$ vanishing for all $a_{1}\geq0$ except possibly at the
poles, which form a subset of $\frac{1}{2}+\mathbb{Z}$. Now let $a_{1}%
\rightarrow-\ell$ with $1\leq\ell\leq a_{3}$. The terms in the left side (the
sum) vanish for $i\leq2\ell-2$ and the right side vanishes for $2a_{1}%
+m+1\leq-1$, that is, $m\leq2\ell-2$. By analytic continuation $\sum
_{i=2\ell-1}^{\min\left(  2a_{2},2a_{3}\right)  }\binom{2a_{1}+2a_{3}%
+1}{2a_{3}-i}\binom{2a_{2}}{i}t_{i}=0$. It remains to change the normalizing
factor in $t_{i}$ to $s\left(  1,2a_{3}+2a_{1}+1,2a_{2}+2a_{1}+1,-2a_{1}%
-1\right)  $. Note that $0\leq c\leq a+b$ implies $2a_{1}+2a_{2}+1\geq0$. The
formal expression (put $i=-2a_{1}-1=2\ell-1$)%
\begin{align*}
&  \frac{s\left(  1,2a_{3}+2a_{1}+1,2a_{2}+2a_{1}+1,-2a_{1}-1\right)
}{s\left(  2a_{1}+2,2a_{3},2a_{2},0\right)  }\\
&  =\frac{\left(  \frac{1}{2}-\kappa-a_{2}\right)  _{\ell-1}\left(  \frac
{1}{2}-\kappa-a_{3}\right)  _{\ell-1}\left(  \kappa+a_{1}+1\right)  _{2\ell
-1}\left(  \frac{3}{2}+a_{1}\right)  _{2\ell-1}}{\left(  \frac{1}{2}%
-a_{2}\right)  _{\ell-1}\left(  \frac{1}{2}-a_{3}\right)  _{\ell-1}\left(
\frac{1}{2}-\kappa-a_{2}-a_{3}\right)  _{2\ell-1}\left(  1-\kappa-a_{2}%
-a_{3}\right)  _{2\ell-1}}%
\end{align*}
has no poles or zeros at integer values of $a_{1}$. Thus the identity remains
valid when multiplied by this ratio.
\end{proof}

The formula for the partial sum in the Theorem was discovered by
experimentation, and recognizing that the factorization of the partial sum
produces linear factors and an irreducible polynomial in $\left(
\kappa+1\right)  \left(  \kappa+n\right)  $. The validity is proved by
induction and a contiguity relation for $F$ (that is, to show $\sum_{i=1}%
^{m}c_{i}=d_{m}$ for sequences $\left\{  c_{i}\right\}  ,\left\{
d_{m}\right\}  $ it suffices to show $d_{0}=0$ and $c_{m}+d_{m-1}-d_{m}=0$ for
$m\geq1$). We use a similar approach to the formula%
\begin{equation}
\left(  8\kappa+2n\right)  \xi_{0}\left(  \left(  q_{1}q_{4}-q_{2}%
q_{3}\right)  \left\langle x\tau\left(  q\right)  ,y\right\rangle
^{2n}\right)  -\xi_{0}\left(  D_{3}\left\langle x\tau\left(  q\right)
,y\right\rangle ^{2n}\right)  =0 \label{d4eqn}%
\end{equation}
for $n=1,2,3,\ldots$. The expression $\left(  8\kappa+2n\right)  \xi
_{0}\left(  \left(  q_{1}q_{4}-q_{2}q_{3}\right)  P_{a,b}^{c}\left(  q\right)
\right)  -\xi_{0}\left(  D_{3}P_{a,b}^{c}\left(  q\right)  \right)  $ with
$a+b=2n$ vanishes when $a$ and $b$ are even, so it suffices to take $a,b,c$
all odd. Furthermore assume $a\geq c$, otherwise use the symmetry $P_{a,b}%
^{c}\left(  q_{4},q_{3},q_{2},q_{1}\right)  =P_{b,a}^{a+b-c}\left(  q\right)
$, which produces the same integral because of the invariance properties of
$q_{1}q_{4}-q_{2}q_{3}$ and $\mu$, and replace $a,b,c$ by $b,a,a+b-c$
respectively. Note that $a+b-c$ is odd and $c\geq a$ implies $a+b-c\leq b$.
Let $a=2a_{1}+2a_{3}+1,b=2a_{2}+1,c=2a_{3}+1$, so $n=a_{1}+a_{2}+a_{3}+1$ and
$a_{1},a_{2},a_{3}\geq0$. The monomials in $P_{a,b}^{c}$ are $q^{\alpha}$ with
$\alpha=\left(  2a_{1}+i,2a_{3}+1-i,2a_{2}+1-i,i\right)  $ and $\left\vert
\alpha\right\vert =2n$. Set%
\begin{gather*}
s\left(  2a_{1},2a_{2}+2,2a_{3}+2,0\right)  t_{i}=\\
\xi_{0}\left(  \left(  \left(  8\kappa+2n\right)  \left(  q_{1}q_{4}%
-q_{2}q_{3}\right)  -D_{3}\right)  q_{1}^{2a_{1}+i}q_{2}^{2a_{3}+1-i}%
q_{3}^{2a_{2}+1-i}q_{4}^{i}\right) \\
=\left(  8\kappa+2n\right)  s\left(  2a_{1}+1+i,2a_{3}+1-i,2a_{2}%
+1-i,i+1\right) \\
-\left(  8\kappa+2n\right)  s\left(  2a_{1}+i,2a_{3}+2-i,2a_{2}+2-i,i\right)
\\
-\left(  2a_{1}+i\right)  s\left(  2a_{1}+i-1,2a_{3}+1-i,2a_{2}+1-i,i+1\right)
\\
-i~s\left(  2a_{1}+i+1,2a_{3}+1-i,2a_{2}+1-i,i-1\right) \\
+\left(  2a_{3}+1-i\right)  s\left(  2a_{1}+i,2a_{3}-i,2a_{2}+2-i,i\right) \\
+\left(  2a_{2}+1-i\right)  s\left(  2a_{1}+i,2a_{3}+2-i,2a_{2}-i,i\right)  .
\end{gather*}
The proof of the following is in Section 3.

\begin{theorem}
\label{big2}Suppose $a_{1}\geq0$ and $m=0,1,2,\ldots$then%
\begin{align*}
\sum_{i=0}^{m}\frac{\left(  -2a_{2}-1\right)  _{i}\left(  -2a_{3}-1\right)
_{i}}{i!\left(  2a_{1}+1\right)  _{i}}t_{i}  &  =\frac{2^{2m+3}\kappa\left(
\kappa+a_{1}\right)  _{m+1}\left(  -2a_{2}\right)  _{m}\left(  -2a_{3}\right)
_{m}\left(  a_{1}+\frac{1}{2}\right)  _{m}}{m!\left(  -2\kappa-2a_{2}%
-2a_{3}-3\right)  _{2m+2}\left(  2a_{1}+1\right)  _{m}}\\
&  \times\left(  4\kappa+3n+2\right)  F\left(  m;\kappa+1,a_{1},a_{2}%
,a_{3}\right)  .
\end{align*}

\end{theorem}

\begin{corollary}
Equation (\ref{d4eqn}) is valid.
\end{corollary}

\begin{proof}
Multiply the formula by $\binom{2a_{1}+2a_{3}+1}{2a_{3}+1}$ and set
$m=\min\left(  2a_{2}+1,2a_{3}+1\right)  $. The factor $\left(  -2a_{2}%
\right)  _{m}\left(  -2a_{3}\right)  _{m}$ vanishes for $m\geq2a_{2}+1$ or
$m\geq2a_{3}+1$.
\end{proof}

In the paper \cite{D3} the proof depended heavily on several integrations by
parts. It may be possible that such a proof exists in this case (given
sufficient ingenuity), but the method of integrating $P_{a,b}^{c}\left(
q\right)  $ by use of $_{4}F_{3}$-series seemed more straightforward.

\section{Hypergeometric series tools}

\subsection{The single sum}

A fundamental transformation for terminating $_{3}F_{2}$ series (for
$n\in\mathbb{N}_{0}$) is%
\[
_{3}F_{2}\left(
\genfrac{}{}{0pt}{}{-n,a,b}{c,d}%
;1\right)  =\frac{\left(  d-b\right)  _{n}}{\left(  d\right)  _{n}}~_{3}%
F_{2}\left(
\genfrac{}{}{0pt}{}{-n,c-a,b}{c,1+b-d-n}%
;1\right)  .
\]
By iterating we obtain two useful transformations:%
\begin{equation}
_{3}F_{2}\left(
\genfrac{}{}{0pt}{}{-n,a,b}{c,d}%
;1\right)  =\frac{\left(  c+d-a-b\right)  _{n}}{\left(  d\right)  _{n}}%
~_{3}F_{2}\left(
\genfrac{}{}{0pt}{}{-n,c-a,c-b}{c,c-a-b+d}%
;1\right)  , \label{chg32a}%
\end{equation}
and (this one provides the sum for the balanced case: $-n+a+b+1=c+d$)%
\begin{align}
_{3}F_{2}\left(
\genfrac{}{}{0pt}{}{-n,a,b}{c,d}%
;1\right)   &  =\left(  -1\right)  ^{n}\frac{\left(  d-a\right)  _{n}\left(
d-b\right)  _{n}}{\left(  c\right)  _{n}\left(  d\right)  _{n}} \label{chg32b}%
\\
&  \times~_{3}F_{2}\left(
\genfrac{}{}{0pt}{}{-n,a+b-n+1-c-d,1-d-n}{a-d+1-n,b-d+1-n}%
;1\right)  .\nonumber
\end{align}
We also need the Whipple transformation (see \cite[p. 56]{B}) for balanced
terminating $_{4}F_{3}$ series: if $n\in\mathbb{N}_{0}$ and $-n+a+b+c+1=c+d+e$
then
\begin{align}
_{4}F_{3}\left(
\genfrac{}{}{0pt}{}{-n,a,b,c}{d,e,f}%
;1\right)   &  =\frac{\left(  e-a\right)  _{n}\left(  f-a\right)  _{n}%
}{\left(  e\right)  _{n}\left(  f\right)  _{n}}\label{Whipple}\\
&  \times~_{4}F_{3}\left(
\genfrac{}{}{0pt}{}{-n,a,d-b,d-c}{d,e-b-c+d,f-b-c+d}%
;1\right)  .\nonumber
\end{align}
With the aid of these formulae we prove the single-sum expression for
$s\left(  \alpha\right)  $. To concentrate on the intermediate steps we change
some variables:%
\begin{align*}
u_{1}  &  =\frac{1}{2}-b_{1}=\frac{1}{2}\left(  1-\alpha_{1}-\alpha
_{4}\right)  ,\\
u_{2}  &  =\frac{1}{2}-b_{0}=\frac{1}{2}\left(  1-\alpha_{2}-\alpha
_{3}\right)  ,\\
m  &  =\frac{1}{2}\left(  \alpha_{3}-\alpha_{4}\right)  ,\\
n  &  =\left\lfloor \frac{\alpha_{4}}{2}\right\rfloor .
\end{align*}
Also write $a$ for $\alpha_{4}$. By the parity condition $m$ is an integer;
further the sets $\left\{  \frac{a}{2},\frac{a-1}{2}\right\}  $ and $\left\{
n,a-n-\frac{1}{2}\right\}  $ are equal. For now we assume $m>-n$ (this is
necessary because factors $\left(  m+1\right)  _{i}$ with $0\leq i\leq n$ will
appear in denominators). The special case $m=-n$, corresponding to $\alpha
_{3}=0,1$ will be handled later. Start with the double sum from equation
(\ref{dblsum}); note $\left(  -\alpha_{4}\right)  _{2i}=2^{2i}\left(
-\frac{\alpha_{4}}{2}\right)  _{i}\left(  \frac{1-\alpha_{4}}{2}\right)
_{i}=2^{2i}\left(  -n\right)  _{i}\left(  n-a+\frac{1}{2}\right)  _{i}$ and
similarly $\left(  -\alpha_{3}\right)  _{2j}=2^{2j}\left(  -m-n\right)
_{j}\left(  -m+n-a+\frac{1}{2}\right)  _{j}$. Further $\frac{1}{2}-b_{3}%
=\frac{1}{2}\left(  1-\alpha_{3}-\alpha_{4}\right)  =\frac{1}{2}-m-a$.%
\begin{align*}
S  &  =\sum_{i=0}^{n}\sum_{j=0}^{m+n}\frac{\left(  -n\right)  _{i}\left(
n-a+\frac{1}{2}\right)  _{i}\left(  -m-n\right)  _{j}\left(  -m+n-a+\frac
{1}{2}\right)  _{j}\left(  \kappa\right)  _{i+j}}{i!j!\left(  u_{1}\right)
_{i}\left(  u_{2}\right)  _{j}\left(  \frac{1}{2}-m-a\right)  _{i+j}}\\
&  =\sum_{i=0}^{n}\frac{\left(  -n\right)  _{i}\left(  n-a+\frac{1}{2}\right)
_{i}\left(  \kappa\right)  _{i}}{i!\left(  u_{1}\right)  _{i}\left(  \frac
{1}{2}-m-a\right)  _{i}}~_{3}F_{2}\left(
\genfrac{}{}{0pt}{}{-m-n,-m+n-a+\frac{1}{2},\kappa+i}{\frac{1}{2}-m-a+i,u_{2}}%
;1\right) \\
&  =f_{1}\sum_{i=0}^{n}\frac{\left(  -n\right)  _{i}\left(  n-a+\frac{1}%
{2}\right)  _{i}\left(  \kappa\right)  _{i}}{i!\left(  u_{1}\right)
_{i}\left(  \frac{1}{2}-m-a\right)  _{i}}~_{3}F_{2}\left(
\genfrac{}{}{0pt}{}{-m-n,\frac{1}{2}-m-a-\kappa,i-n}{\frac{1}{2}%
-m-a+i,-n-\kappa+u_{2}}%
;1\right) \\
&  =f_{1}\sum_{i=0}^{n}\sum_{j=0}^{n-i}\frac{\left(  -n\right)  _{i+j}\left(
n-a+\frac{1}{2}\right)  _{i}\left(  -m-n\right)  _{j}\left(  \frac{1}%
{2}-m-a-\kappa\right)  _{j}\left(  \kappa\right)  _{i}}{i!j!\left(
u_{1}\right)  _{i}\left(  -n-\kappa+u_{2}\right)  _{j}\left(  \frac{1}%
{2}-m-a\right)  _{i+j}}\\
&  =f_{1}S^{\prime},~f_{1}:=\frac{\left(  -\kappa+u_{2}-n\right)  _{n+m}%
}{\left(  u_{2}\right)  _{n+m}}.
\end{align*}
The transformation (\ref{chg32a}) is used on the inner sum. Note that the
range of summation is now changed to $\left\{  \left(  i,j\right)
:i\geq0,j\geq0,i+j\leq n\right\}  $, and $m$ can be considered a formal
parameter. Reversing the order of summation and using equation (\ref{chg32b})
we obtain%
\begin{align*}
S^{\prime}  &  =\sum_{j=0}^{n}\frac{\left(  -n\right)  _{j}\left(
-m-n\right)  _{j}\left(  \frac{1}{2}-m-a-\kappa\right)  _{j}}{j!\left(
-n-\kappa+u_{2}\right)  _{j}\left(  \frac{1}{2}-m-a\right)  _{j}}~_{3}%
F_{2}\left(
\genfrac{}{}{0pt}{}{j-n,n-a+\frac{1}{2},\kappa}{u_{1},\frac{1}{2}-m-a+j}%
;1\right) \\
&  =\sum_{j=0}^{n}\frac{\left(  -n\right)  _{j}\left(  -m-n\right)
_{j}\left(  \frac{1}{2}-m-a-\kappa\right)  _{j}}{j!\left(  -n-\kappa
+u_{2}\right)  _{j}\left(  \frac{1}{2}-m-a\right)  _{j}}\\
&  \times\left(  -1\right)  ^{n-j}\frac{\left(  -m+j-n\right)  _{n-j}\left(
\frac{1}{2}-m-a-\kappa+j\right)  _{n-j}}{\left(  u_{1}\right)  _{n-j}\left(
\frac{1}{2}-m-a+j\right)  _{n-j}}\\
&  \times~_{3}F_{2}\left(
\genfrac{}{}{0pt}{}{j-n,-u_{1}+\kappa+m+1,\frac{1}{2}+m+a-n}{\frac{1}%
{2}+\kappa+m+a-n,m+1}%
;1\right)  ,
\end{align*}
and reversing the order of summation again we obtain%
\begin{align*}
S^{\prime}  &  =\left(  -1\right)  ^{n}\frac{\left(  -m-n\right)  _{n}\left(
\frac{1}{2}-m-a-\kappa\right)  _{n}}{\left(  \frac{1}{2}-m-a\right)
_{n}\left(  u_{1}\right)  _{n}}\\
&  \times\sum_{i=0}^{n}\sum_{j=0}^{n-i}\frac{\left(  -n\right)  _{i}\left(
-u_{1}+\kappa+m+1\right)  _{i}\left(  \frac{1}{2}+m+a-n\right)  _{i}\left(
i-n\right)  _{j}\left(  1-n-u_{1}\right)  _{j}}{i!j!\left(  \frac{1}{2}%
+\kappa+m+a-n\right)  _{i}\left(  m+1\right)  _{i}\left(  -n-\kappa
+u_{2}\right)  _{j}}\\
&  =f_{2}\sum_{i=0}^{n}\frac{\left(  -n\right)  _{i}\left(  -u_{1}%
+\kappa+m+1\right)  _{i}\left(  \frac{1}{2}+m+a-n\right)  _{i}\left(
u_{1}+u_{2}-\kappa-1\right)  _{n-i}}{i!\left(  \frac{1}{2}+\kappa
+m+a-n\right)  _{i}\left(  m+1\right)  _{i}\left(  -n-\kappa+u_{2}\right)
_{n-i}}\\
&  =f_{2}S^{\prime\prime},~f_{2}:=\frac{\left(  m+1\right)  _{n}\left(
\frac{1}{2}-m-a-\kappa\right)  _{n}}{\left(  \frac{1}{2}-m-a\right)
_{n}\left(  u_{1}\right)  _{n}}.
\end{align*}

At the second last equation the $j$-sum is done with the Chu-Vandermonde sum
($_{2}F_{1}\left(
\genfrac{}{}{0pt}{}{-n,b}{c}%
;1\right)  =\frac{\left(  c-b\right)  _{n}}{\left(  c\right)  _{n}}$). Replace
$i$ by $n-i$ in the sum $S^{\prime\prime}$, and use the relations $\left(
b\right)  _{n-i}=\left(  -1\right)  ^{i}\left(  b\right)  _{n}/\left(
1-b-n\right)  _{i}$ and $\left(  -n\right)  _{n-i}/\left(  n-i\right)
!=\left(  -1\right)  ^{n}\left(  -n\right)  _{i}/i!$ to obtain%
\begin{align*}
S^{\prime\prime}  &  =\left(  -1\right)  ^{n}\frac{\left(  -u_{1}%
+\kappa+m+1\right)  _{n}\left(  \frac{1}{2}+m+a-n\right)  _{n}}{\left(
\frac{1}{2}+\kappa+m+a-n\right)  _{n}\left(  m+1\right)  _{n}}\\
&  \times~_{4}F_{3}\left(
\genfrac{}{}{0pt}{}{-n,u_{1}+u_{2}-\kappa-1,\frac{1}{2}-m-a-\kappa
,-m-n}{\frac{1}{2}-m-a,u_{1}-\kappa-m-n,-n-\kappa+u_{2}}%
;1\right) \\
&  =f_{3}~_{4}F_{3}\left(
\genfrac{}{}{0pt}{}{-n,u_{1}+u_{2}-\kappa-1,n-a+\frac{1}{2},\kappa}{\frac
{1}{2}-m-a,u_{2}+m,u_{1}}%
;1\right)  ,\\
f_{3}  &  =\left(  -1\right)  ^{n}\frac{\left(  -u_{1}+\kappa+m+1\right)
_{n}\left(  \frac{1}{2}-m-a\right)  _{n}\left(  1-u_{1}-n\right)  _{n}\left(
1-u_{2}-m-n\right)  _{n}}{\left(  \frac{1}{2}-\kappa-m-a\right)  _{n}\left(
m+1\right)  _{n}\left(  -n-\kappa+u_{2}\right)  _{n}\left(  -n+u_{1}%
-\kappa-m\right)  _{n}}.
\end{align*}
The Whipple transformation is valid since the first $_{4}F_{3}$-series is
balanced (as is the second, of course). We combine the factors
\begin{align*}
f_{1}f_{2}f_{3}  &  =\frac{\left(  -\kappa+u_{2}-n\right)  _{n+m}}{\left(
u_{2}\right)  _{n+m}}\frac{\left(  m+1\right)  _{n}\left(  \frac{1}%
{2}-m-a-\kappa\right)  _{n}}{\left(  \frac{1}{2}-m-a\right)  _{n}\left(
u_{1}\right)  _{n}}\\
&  \times\frac{\left(  -u_{1}+\kappa+m+1\right)  _{n}\left(  \frac{1}%
{2}-m-a\right)  _{n}\left(  u_{1}\right)  _{n}\left(  u_{2}+m\right)  _{n}%
}{\left(  \frac{1}{2}-\kappa-m-a\right)  _{n}\left(  m+1\right)  _{n}\left(
-n-\kappa+u_{2}\right)  _{n}\left(  -u_{1}+\kappa+m+1\right)  _{n}}\\
&  =\frac{\left(  u_{2}-\kappa\right)  _{m}}{\left(  u_{2}\right)  _{m}}.
\end{align*}
Thus%
\[
S=\frac{\left(  u_{2}-\kappa\right)  _{m}}{\left(  u_{2}\right)  _{m}}%
~_{4}F_{3}\left(
\genfrac{}{}{0pt}{}{-n,u_{1}+u_{2}-\kappa-1,n-a+\frac{1}{2},\kappa}{\frac
{1}{2}-m-a,u_{2}+m,u_{1}}%
;1\right)  .
\]
It is possible that $0>m>-n$ (that is, $\alpha_{3}<\alpha_{4}$) in which case
the equation $\left(  t\right)  _{m}=\frac{1}{\left(  t+m\right)  _{-m}}%
=\frac{\left(  -1\right)  ^{m}}{\left(  1-t\right)  _{-m}}$ applies, and so
$\dfrac{\left(  u_{2}-\kappa\right)  _{m}}{\left(  u_{2}\right)  _{m}}%
=\dfrac{\left(  1-u_{2}\right)  _{-m}}{\left(  1+\kappa-u_{2}\right)  _{-m}}$.
Finally for the special case $m=-n$ the double sum $S$ reduces to the single
sum%
\begin{align*}
S  &  =\sum_{i=0}^{n}\frac{\left(  -n\right)  _{i}\left(  n-a+\frac{1}%
{2}\right)  _{i}\left(  \kappa\right)  _{i}}{i!\left(  u_{1}\right)
_{i}\left(  \frac{1}{2}-m-a\right)  _{i}}=\sum_{i=0}^{n}\frac{\left(
-n\right)  _{i}\left(  \kappa\right)  _{i}}{i!\left(  u_{1}\right)  _{i}}\\
&  =\frac{\left(  u_{1}-\kappa\right)  _{n}}{\left(  u_{1}\right)  _{n}},
\end{align*}
while the $_{4}F_{3}$ series reduces to a balanced $_{3}F_{2}$ series, so
that
\[
_{3}F_{2}\left(
\genfrac{}{}{0pt}{}{-n,u_{1}+u_{2}-\kappa-1,\kappa}{u_{2}-n,u_{1}}%
\right)  =\left(  -1\right)  ^{n}\frac{\left(  1+\kappa-u_{2}\right)
_{n}\left(  u_{1}-\kappa\right)  _{n}}{\left(  u_{2}-n\right)  _{n}\left(
u_{1}\right)  _{n}},
\]
(by use of formula (\ref{chg32b})) and the result multiplied by $\dfrac
{\left(  1-u_{2}\right)  _{-m}}{\left(  1+\kappa-u_{2}\right)  _{-m}}$
produces $\dfrac{\left(  u_{1}-\kappa\right)  _{n}}{\left(  u_{1}\right)
_{n}}$, which equals $S$.

This finishes the argument, and now we can show the three-term recurrence
applies to $s\left(  \alpha\right)  .$ Returning to the original variables, we
note $b_{0}-m=\frac{1}{2}\left(  \alpha_{2}+\alpha_{3}-\alpha_{3}+\alpha
_{4}\right)  =b_{2}$, thus%
\begin{align*}
\dfrac{\left(  u_{2}-\kappa\right)  _{m}}{\left(  u_{2}\right)  _{m}}  &
=\frac{\left(  \frac{1}{2}+b_{0}+\kappa-m\right)  _{m}\left(  \kappa+\frac
{1}{2}\right)  _{b_{2}}\left(  \frac{1}{2}\right)  _{b_{2}}}{\left(  \frac
{1}{2}+b_{0}-m\right)  _{m}\left(  \frac{1}{2}\right)  _{b_{2}}\left(
\kappa+\frac{1}{2}\right)  _{b_{2}}}\\
&  =\frac{\left(  \kappa+\frac{1}{2}\right)  _{b_{0}}\left(  \frac{1}%
{2}\right)  _{b_{2}}}{\left(  \frac{1}{2}\right)  _{b_{0}}\left(  \kappa
+\frac{1}{2}\right)  _{b_{2}}}.
\end{align*}
This proves formula (\ref{sngl}).

\subsection{Three term recurrence and symmetries}

\begin{proposition}
For $n\in\mathbb{N}_{0}$
\begin{align*}
&  n\left(  n+2v_{1}\right)  \left(  n+\kappa+v_{1}-\frac{1}{2}\right)
\left(  \kappa+\frac{1}{2}+v_{2}+v_{3}-n\right)  F\left(  n-1;\kappa
,v_{1},v_{2},v_{3}\right) \\
&  +\left\{  n\left(  n+2v_{1}\right)  \left(  n-v_{2}-v_{3}-\frac{1}%
{2}\right)  +\left(  n+\frac{1}{2}+v_{1}\right)  \left(  n-2v_{2}\right)
\left(  n-2v_{3}\right)  \right\} \\
&  \times\left(  n+v_{1}-\frac{1}{2}\right)  F\left(  n;\kappa,v_{1}%
,v_{2},v_{3}\right) \\
&  =\left(  n+v_{1}-\frac{1}{2}\right)  \left(  n+v_{1}+\frac{1}{2}\right)
\left(  n-2v_{2}\right)  \left(  n-2v_{3}\right)  F\left(  n+1;\kappa
,v_{1},v_{2},v_{3}\right) \\
&  .
\end{align*}

\end{proposition}

\begin{proof}
Expand the left side as a series in $g_{i}:=\left(  \kappa\right)  _{i}\left(
-\kappa-v_{1}-v_{2}-v_{3}\right)  _{i}$ for $0\leq i\leq\left\lfloor
\frac{n+1}{2}\right\rfloor $. Start with the equation%
\begin{gather*}
\left(  n+\kappa+v_{1}-\frac{1}{2}\right)  \left(  \kappa+\frac{1}{2}%
+v_{2}+v_{3}-n\right)  =\\
\left(  \frac{1}{2}-v_{1}-n+i\right)  \left(  n-v_{2}-v_{3}-\frac{1}%
{2}+i\right)  -\left(  \kappa+i\right)  \left(  -\kappa-v_{1}-v_{2}%
-v_{3}+i\right)  .
\end{gather*}

The outline of the calculation is this: write $F\left(  n-1\right)
=\sum_{i\geq0}a_{i}g_{i}$ and $F\left(  n\right)  =\sum_{i\geq0}a_{i}^{\prime
}g_{i}$ (suppressing the other arguments); note that $g_{i}\left(
\kappa+i\right)  \left(  -\kappa-v_{1}-v_{2}-v_{3}+i\right)  =g_{i+1}$ then
collect term-by-term in
\begin{align*}
&  \left(  n+2v_{1}\right)  \{n\left(  \frac{1}{2}-v_{1}-n+i\right)  \left(
n-v_{2}-v_{3}-\frac{1}{2}+i\right)  \sum_{i\geq0}a_{i}g_{i}\\
&  -n\sum_{i\geq1}a_{i-1}g_{i}+n\left(  n-v_{2}-v_{3}-\frac{1}{2}\right)
\left(  n+v_{1}-\frac{1}{2}\right)  \sum_{i\geq0}a_{i}^{\prime}g_{i}\}\\
&  +\left(  n+\frac{1}{2}+v_{1}\right)  \left(  n+v_{1}-\frac{1}{2}\right)
\left(  n-2v_{2}\right)  \left(  n-2v_{3}\right)  \sum_{i\geq0}a_{i}^{\prime
}g_{i}.
\end{align*}
Write this sum as $\sum_{i\geq0}\frac{\left(  -\frac{n}{2}\right)  _{i}\left(
-\frac{1+n}{2}\right)  _{i}}{i!\left(  \frac{1}{2}-n-v_{1}\right)  _{i}\left(
\frac{1}{2}-v_{2}\right)  _{i}\left(  \frac{1}{2}-v_{3}\right)  _{i}}%
c_{i}g_{i}$, then (in corresponding order)%
\begin{align*}
c_{i}  &  =\frac{n+2v_{1}}{n+1}\left(  n+v_{1}-\frac{1}{2}\right)  \{-\left(
n-2i\right)  \left(  n+1-2i\right)  \left(  n-v_{2}-v_{3}-\frac{1}{2}+i\right)
\\
&  +n\left(  n+1-2i\right)  \left(  n-v_{2}-v_{3}-\frac{1}{2}\right)
+i\left(  2v_{2}+1-2i\right)  \left(  2v_{3}+1-i\right)  \}\\
&  +\frac{n+1-2i}{n+1}\left(  n+v_{1}-\frac{1}{2}\right)  \left(
n-2v_{2}\right)  \left(  n-2v_{3}\right) \\
&  =\frac{\left(  n-2v_{2}\right)  \left(  n-2v_{3}\right)  }{n+1}\left(
n+v_{1}-\frac{1}{2}\right)  \left\{  i\left(  n+2v_{1}\right)  +\left(
n+v_{1}+\frac{1}{2}\right)  \left(  n+1-2i\right)  \right\} \\
&  =\left(  n-2v_{2}\right)  \left(  n-2v_{3}\right)  \left(  n+v_{1}-\frac
{1}{2}\right)  \left(  n+v_{1}+\frac{1}{2}+i\right)  .
\end{align*}
Finally $\left(  n+v_{1}+\frac{1}{2}+i\right)  /\left(  \frac{1}{2}%
-n-v_{1}\right)  _{i}=\left(  n+v_{1}+\frac{1}{2}\right)  /\left(  -\frac
{1}{2}-n-v_{1}\right)  _{i}$, \ which proves the identity.
\end{proof}

\begin{corollary}
The recurrence (\ref{recurs}) for $s^{\prime}\left(  \alpha\right)  $ is valid.
\end{corollary}

\begin{proof}
Suppose $\alpha_{1}\equiv\alpha_{2}\equiv\alpha_{3}\equiv\alpha_{4}%
\operatorname{mod}2$ then set $c^{\prime}\left(  \alpha\right)  =\prod
_{i=1}^{3}\frac{\left(  \frac{1}{2}\right)  _{b_{i}}}{\left(  \kappa+\frac
{1}{2}\right)  _{b_{i}}}$ where $b_{i}=\frac{1}{2}\left(  \alpha_{i}%
+\alpha_{4}\right)  ,1\leq i\leq3$. Then $s^{\prime}\left(  \alpha\right)
/c^{\prime}\left(  \alpha\right)  $ satisfies the recurrence for $F$ in the
Proposition, with $n=\alpha_{4},v_{1}=\frac{1}{2}\left(  \alpha_{1}-\alpha
_{4}\right)  ,v_{2}=b_{2},v_{3}=b_{3}$. Multiply the recurrence by $c^{\prime
}\left(  \alpha\right)  $ then the coefficient of $F\left(  n-1\right)  $ is
multiplied by%
\[
\dfrac{c^{\prime}\left(  \alpha\right)  }{c^{\prime}\left(  \alpha
_{1}-1,\alpha_{2}+1,\alpha_{3}+1,\alpha_{4}-1\right)  }=\frac{\alpha
_{1}+\alpha_{4}-1}{2\kappa+\alpha_{1}+\alpha_{4}-1}%
\]
and the coefficient of $F\left(  n+1\right)  $ is multiplied by%
\[
\dfrac{c^{\prime}\left(  \alpha\right)  }{c^{\prime}\left(  \alpha
_{1}+1,\alpha_{2}-1,\alpha_{3}-1,\alpha_{4}+1\right)  }=\frac{2\kappa
+\alpha_{1}+\alpha_{4}+1}{\alpha_{1}+\alpha_{4}+1}.
\]
Divide out the common factor $\frac{1}{2}\left(  \alpha_{1}+\alpha
_{4}-1\right)  =\left(  n+v_{1}-\frac{1}{2}\right)  $ to obtain (\ref{recurs}).
\end{proof}

There is another, perhaps unexpected, symmetry:

\begin{proposition}
$s^{\prime}$ is completely symmetric in its arguments.
\end{proposition}

\begin{proof}
Already formula \ref{sngl} shows the symmetry in $\left(  \alpha_{1}%
,\alpha_{2},\alpha_{3}\right)  $. One can argue from the $\alpha
_{2}\leftrightarrow\alpha_{3}$ and $\alpha_{1}\leftrightarrow\alpha_{4}$
invariance together with the $\alpha\leftrightarrow\left(  \alpha_{2}%
,\alpha_{1},\alpha_{4},\alpha_{3}\right)  $ invariance of $s$. The Whipple
transformation (\ref{Whipple}) gives a direct proof. It suffices to show
$s^{\prime}\left(  \alpha_{1},\alpha_{2},\alpha_{3},\alpha_{4}\right)
=s^{\prime}\left(  \alpha_{4},\alpha_{2},\alpha_{3},\alpha_{1}\right)  $. If
all the $\alpha_{i}$'s are even take $n=\frac{\alpha_{4}}{2},a=\kappa
,d=\frac{1}{2}-b_{1}$, then%
\begin{align*}
s^{\prime}\left(  \alpha_{1},\alpha_{2},\alpha_{3},\alpha_{4}\right)   &
=\frac{\left(  \frac{1}{2}\right)  _{\alpha_{2}/2}\left(  \frac{1}{2}\right)
_{\alpha_{3}/2}}{\left(  \kappa+\frac{1}{2}\right)  _{\alpha_{2}/2}\left(
\kappa+\frac{1}{2}\right)  _{\alpha_{3}/2}}\times\\
&  _{4}F_{3}\left(
\genfrac{}{}{0pt}{}{-\frac{\alpha_{4}}{2},\kappa,\frac{-\alpha_{1}}{2}%
,\kappa+\frac{1}{2}\left(  1+\alpha_{2}+\alpha_{3}\right)  }{\frac{1}{2}%
-b_{1},\kappa+\frac{1}{2}\left(  1+\alpha_{3}\right)  ,\kappa+\frac{1}%
{2}\left(  1+\alpha_{2}\right)  }%
;1\right)  ,
\end{align*}
and this expression is symmetric in $\left(  \alpha_{1},\alpha_{4}\right)  $.
A similar argument works when $\alpha_{4}$ is odd:%
\begin{align*}
s^{\prime}\left(  \alpha_{1},\alpha_{2},\alpha_{3},\alpha_{4}\right)   &
=\frac{\left(  \frac{1}{2}\right)  _{\left(  \alpha_{2}+1\right)  /2}\left(
\frac{1}{2}\right)  _{\left(  \alpha_{3}+1\right)  /2}}{\left(  \kappa
+\frac{1}{2}\right)  _{\left(  \alpha_{2}+1\right)  /2}\left(  \kappa+\frac
{1}{2}\right)  _{\left(  \alpha_{3}+1\right)  /2}}\times\\
&  _{4}F_{3}\left(
\genfrac{}{}{0pt}{}{\frac{1-\alpha_{4}}{2},\kappa,\frac{1-\alpha_{1}}%
{2},\kappa+\frac{1}{2}\left(  1+\alpha_{2}+\alpha_{3}\right)  }{\frac{1}%
{2}-b_{1},\kappa+1+\frac{\alpha_{3}}{2},\kappa+1+\frac{\alpha_{2}}{2}}%
;1\right)  .
\end{align*}

\end{proof}

\subsection{\label{sctcont}Contiguity Relations}

We start by writing the expressions in Theorems \ref{big1} and \ref{big2} in
terms of the function $F$ (see definition \ref{defF}). Suppose $\alpha
_{1}\equiv\alpha_{2}\equiv\alpha_{3}\equiv\alpha_{4}\operatorname{mod}2$ then
set%
\[
c\left(  \alpha\right)  =\frac{\left(  2\kappa\right)  _{\alpha_{1}+\alpha
_{4}}\left(  2\kappa\right)  _{\alpha_{2}+\alpha_{3}}}{\left(  4\kappa\right)
_{\left\vert \alpha\right\vert }}\prod_{i=1}^{3}\frac{\left(  \frac{1}%
{2}\right)  _{b_{i}}}{\left(  \kappa+\frac{1}{2}\right)  _{b_{i}}}%
\]
where $b_{i}=\frac{1}{2}\left(  \alpha_{i}+\alpha_{4}\right)  ,1\leq i\leq3$;
thus $s\left(  \alpha\right)  =c\left(  \alpha\right)  F\left(  \alpha
_{4};\kappa,\frac{1}{2}\left(  \alpha_{1}-\alpha_{4}\right)  ,b_{2}%
,b_{3}\right)  $. For Theorem \ref{big1} we set $i=m\geq1$ and extract the
common factor%
\begin{align*}
&  C_{m}=\frac{c\left(  2a_{1}+m,2a_{3}-m,2a_{3}-m,m\right)  \left(
-2a_{3}\right)  _{m}\left(  -2a_{2}\right)  _{m}}{s\left(  2a_{1}%
+2,2a_{3},2a_{2},0\right)  m!\left(  2a_{1}+2\right)  _{m}\left(
4\kappa+n\right)  \left(  2\kappa+n\right)  }\\
&  =2^{2m-1}\frac{\left(  4\kappa+2n+1\right)  \left(  \kappa+a_{1}+1\right)
_{m-1}\left(  a_{1}+\frac{3}{2}\right)  _{m-1}\left(  -2a_{3}\right)
_{m}\left(  -2a_{2}\right)  _{m}}{\left(  4\kappa+n\right)  \left(
1-2\kappa-2a_{2}-2a_{3}\right)  _{2m}m!\left(  2a_{1}+2\right)  _{m}}.
\end{align*}
Recall $n=\sum_{i=1}^{3}a_{i}$. When $i=0$ we have
\begin{align*}
t_{0}  &  =2\left(  4\kappa+2n+1\right)  -2\frac{\left(  2a_{1}+1\right)
\left(  \kappa+a_{1}\right)  c\left(  2a_{1},2a_{3},2a_{3},0\right)  }{\left(
2\kappa+n\right)  s\left(  2a_{1}+2,2a_{3},2a_{2},0\right)  }\\
&  =0.
\end{align*}
The plan is to use induction by showing $\frac{\left(  -2a_{3}\right)
_{m}\left(  -2a_{2}\right)  _{m}}{m!\left(  2a_{1}+2\right)  _{m}}%
t_{m}+d_{m-1}-d_{m}=0$, where $d_{m}$ denotes the claimed expression for the
sum. The following multipliers are needed in the calculation of $t_{m}$:%
\begin{align*}
\frac{c\left(  2a_{1}+m+2,2a_{3}-m,2a_{2}-m,m\right)  }{c\left(
2a_{1}+m,2a_{3}-m,2a_{2}-m,m\right)  }  &  =\frac{\left(  \kappa
+a_{1}+m\right)  \left(  2a_{1}+2m+1\right)  }{\left(  4\kappa+2n+1\right)
\left(  2\kappa+n\right)  },\\
\frac{c\left(  2a_{1}+m+2,2a_{3}-m,2a_{2}-m,m-2\right)  }{c\left(
2a_{1}+m,2a_{3}-m,2a_{2}-m,m\right)  }  &  =\frac{\left(  2\kappa
+2a_{2}-1\right)  \left(  2\kappa+2a_{3}-1\right)  }{\left(  2a_{2}-1\right)
\left(  2a_{3}-1\right)  },\\
\frac{c\left(  2a_{1}+m+1,2a_{3}-m-1,2a_{2}-m+1,m-1\right)  }{c\left(
2a_{1}+m,2a_{3}-m,2a_{2}-m,m\right)  }  &  =\frac{2\kappa+2a_{3}-1}{2a_{3}%
-1},\\
\frac{c\left(  2a_{1}+m+1,2a_{3}-m+1,2a_{2}-m-1,m-1\right)  }{c\left(
2a_{1}+m,2a_{3}-m,2a_{2}-m,m\right)  }  &  =\frac{2\kappa+2a_{2}-1}{2a_{2}-1}.
\end{align*}
Write $d_{m}=A_{m}F\left(  m-1;\kappa+1,a_{1}+1,a_{2}-1,a_{3}-1\right)  $
then
\begin{align*}
\frac{A_{m}}{C_{m}}  &  =4m\kappa\left(  \kappa+a_{1}+m\right)  \frac{\left(
2a_{2}-m\right)  \left(  2a_{3}-m\right)  }{\left(  2a_{2}-1\right)  \left(
2a_{3}-1\right)  },\\
\frac{A_{m-1}}{A_{m}}  &  =\frac{2\left(  m-1\right)  \left(  2a_{1}%
+m+1\right)  \left(  \kappa+a_{2}+a_{3}-m\right)  \left(  \kappa+a_{2}%
+a_{3}+\frac{1}{2}-m\right)  }{\left(  \kappa+a_{1}+m\right)  \left(
2a_{1}+2m-1\right)  \left(  2a_{2}-m\right)  \left(  2a_{3}-m\right)  }.
\end{align*}
The expression $\left(  \frac{\left(  -2a_{3}\right)  _{m}\left(
-2a_{2}\right)  _{m}}{m!\left(  2a_{1}+2\right)  _{m}}t_{m}+d_{m-1}%
-d_{m}\right)  /C_{m}$ with each $a_{i}$ is replaced by $v_{i}$ (formal
parameters instead of integers) becomes the left side of the equation in the following.

\begin{theorem}
Suppose $v_{1},v_{2},v_{3}\notin-\frac{1}{2}+\mathbb{N}_{0}$, ($v_{0}%
=v_{1}+v_{2}+v_{3}$) and $m=1,2,3,\ldots$then%
\begin{align*}
&  2\left(  2v_{1}+2m+1\right)  \left(  \kappa+v_{1}+m\right)  \left(
4\kappa+v_{0}\right)  F\left(  m;\kappa,v_{1}+1,v_{2},v_{3}\right) \\
&  -\left\{  2\left(  \kappa+v_{1}+m\right)  \left(  4\kappa+v_{0}\right)
+m\left(  3\kappa+v_{0}\right)  \right\}  \left(  2v_{1}+m+1\right)  F\left(
m;\kappa,v_{1},v_{2},v_{3}\right) \\
&  -m\left(  m-1\right)  \frac{\left(  2v_{2}+2\kappa-1\right)  \left(
2v_{3}+2\kappa-1\right)  }{\left(  2v_{2}-1\right)  \left(  2v_{3}-1\right)
}\left(  3\kappa+v_{0}\right)  F\left(  m-2;\kappa,v_{1}+2,v_{2}%
-1,v_{3}-1\right) \\
&  -m\kappa\frac{\left(  2\kappa+2v_{3}-1\right)  }{2v_{3}-1}\left(
2v_{3}-m\right)  F\left(  m-1;\kappa,v_{1}+1,v_{2},v_{3}-1\right) \\
&  -m\kappa\frac{\left(  2\kappa+2v_{2}-1\right)  }{2v_{2}-1}\left(
2v_{2}-m\right)  F\left(  m-1;\kappa,v_{1}+1,v_{2}-1,v_{3}\right) \\
&  +m\left(  m-1\right)  \frac{8\kappa\left(  2v_{1}+m+1\right)  }{\left(
2v_{2}-1\right)  \left(  2v_{3}-1\right)  \left(  2v_{1}+2m-1\right)  }\left(
\kappa+v_{2}+v_{3}-m\right)  \left(  \kappa+v_{2}+v_{3}+\frac{1}{2}-m\right)
\\
&  \times F\left(  m-2;\kappa+1,v_{1}+1,v_{2}-1,v_{3}-1\right) \\
&  -4m\kappa\frac{\left(  2v_{2}-m\right)  \left(  2v_{3}-m\right)  }{\left(
2v_{2}-1\right)  \left(  2v_{3}-1\right)  }\left(  \kappa+v_{1}+m\right)
F\left(  m-1;\kappa+1,v_{1}+1,v_{2}-1,v_{3}-1\right) \\
&  =0.
\end{align*}

\end{theorem}

\begin{proof}
We proceed by considering the left side as a polynomial in $\kappa$ of the
form%
\[
\sum_{i\geq0}\left(  c_{i}+d_{i}\kappa\right)  \left(  \kappa\right)
_{i}\left(  -\kappa-v_{1}-v_{2}-v_{3}\right)  _{i}.
\]
The coefficients are rational functions of $v_{1},v_{2},v_{3},m$. Write the
identity in abbreviated form as $\sum_{n=1}^{7}\mu_{n}F\left(  \gamma
_{n}\right)  =0$ using the same order as above (for example, $\gamma
_{2}=\left(  m;\kappa,v_{1},v_{2},v_{3}\right)  $ and $\mu_{5}=-m\kappa
\frac{\left(  2\kappa+2v_{2}-1\right)  }{2v_{2}-1}\left(  2v_{2}-m\right)  $).
Denote term $\#i$ of $F\left(  \gamma_{n}\right)  $ by $t_{i}\left(
\gamma_{n}\right)  $, (where $t_{0}\left(  \gamma_{n}\right)  =1$ for all $n$)
so that
\[
t_{i}\left(  \gamma_{3}\right)  =\frac{\left(  1-\frac{m}{2}\right)
_{i}\left(  \frac{3-m}{2}\right)  _{i}\left(  \kappa\right)  _{i}\left(
-\kappa-v_{1}-v_{2}-v_{3}\right)  _{i}}{\left(  \frac{1}{2}-v_{1}-m\right)
_{i}\left(  \frac{3}{2}-v_{2}\right)  _{i}\left(  \frac{3}{2}-v_{3}\right)
_{i}~i!},
\]
for example; with the understanding that $F\left(  m\right)  =0$ for $m<0$.
The verification is carried out by symbolic computation. The underlying idea
is to set up the equations as rational functions of all of the variables. For
example, $\left(  a\right)  _{i}$ is not rational in $i$, but for a fixed
integer $k$ the expression $\left(  a\right)  _{i}/\left(  a+k\right)  _{i}$
is rational in $a$ and $i$ because of the identity $\left(  a\right)
_{i}/\left(  a+k\right)  _{i}=\left(  a\right)  _{k}/\left(  a+i\right)  _{k}%
$. By use of symbolic calculation we evaluate%
\begin{align*}
p_{i}\left(  \kappa\right)   &  =\mu_{1}\frac{t_{i+1}\left(  \gamma
_{1}\right)  }{t_{i}\left(  \gamma_{2}\right)  }+\mu_{2}+\sum_{n=3}^{7}\mu
_{n}\frac{t_{i}\left(  \gamma_{n}\right)  }{t_{i}\left(  \gamma_{2}\right)
},\\
p_{0}\left(  \kappa\right)   &  =\mu_{1}\left(  1+t_{1}\left(  \gamma
_{1}\right)  \right)  +\mu_{2}+\sum_{n=3}^{7}\mu_{n},
\end{align*}
for $1\leq i\leq\left\lfloor \frac{m}{2}\right\rfloor $. Indeed $p_{i}\left(
\kappa\right)  $ is a polynomial in $\kappa$ of degree 4 (note that
$\frac{t_{i}\left(  \gamma_{1}\right)  }{t_{i}\left(  \gamma_{2}\right)  }$ is
not polynomial in $\kappa$ for $i\geq1$, and the factor $\kappa$ in $\mu_{6}$
and $\mu_{7}$ is necessary since $\frac{t_{i}\left(  \gamma_{n}\right)
}{t_{i}\left(  \gamma_{2}\right)  }$ has $\kappa$ in the denominator for
$n=6,7$), further $p_{i}\left(  \kappa\right)  $ is rational in all the
variables ($v_{1},v_{2},v_{3},\kappa,i,m$). Although $i$ is an indexing
variable the difference between the respective parameters in $\left\{
\gamma_{n}\right\}  $ is one of \thinspace$0,\pm1,\pm2$ and thus $\frac
{t_{i}\left(  \gamma_{n}\right)  }{t_{i}\left(  \gamma_{2}\right)  }$ is
rational in $i$ for each $n$. By computer algebra we find the coefficients
$\left\{  c_{i,k}\right\}  _{k=0}^{4}$, which are functions of $i,m,v_{1}%
,v_{2},v_{3}$, so that%
\[
p_{i}\left(  \kappa\right)  =c_{i,0}+c_{i,1}\kappa+\left(  c_{i,2}%
+c_{i,3}\kappa\right)  g_{i,1}\left(  \kappa\right)  +c_{i,4}g_{i,2}\left(
\kappa\right)  ,
\]
where $g_{i,k}\left(  \kappa\right)  :=\left(  \kappa+i\right)  _{k}\left(
-\kappa-v_{0}+i\right)  _{k}$ for $k\geq0$, thus $g_{0,i}\left(
\kappa\right)  g_{i,k}\left(  \kappa\right)  =g_{0,i+k}\left(  \kappa\right)
$. For example
\[
c_{i,4}=\frac{4\left(  2i-m\right)  \left(  2i+1-m\right)  }{\left(
i+1\right)  \left(  \frac{1}{2}-v_{2}+i\right)  \left(  \frac{1}{2}%
-v_{3}+i\right)  }%
\]
for $0\leq i\leq\left\lfloor \frac{m}{2}\right\rfloor $, but the other
coefficients are more complicated. Because the degree of $\kappa$ in $\mu_{1}$
is 2 the value of $c_{0,4}$ agrees with the generic $c_{i,4}$ with $i=0$. For
$i\geq0$ let
\[
r_{i}=\frac{t_{i}\left(  \gamma_{2}\right)  }{g_{0,i}\left(  \kappa\right)
},
\]
so that $r_{i}$ is independent of $\kappa$ (thus $r_{i}=\dfrac{\left(
-\frac{m}{2}\right)  _{i}\left(  \frac{1-m}{2}\right)  _{i}}{i!\left(
\frac{1}{2}-v_{1}-m\right)  _{i}\left(  \frac{1}{2}-v_{2}\right)  _{i}\left(
\frac{1}{2}-v_{3}\right)  _{i}}$). The left side of the identity equals%
\begin{align*}
&  \sum_{i=0}^{\left\lfloor m/2\right\rfloor }p_{i}\left(  \kappa\right)
g_{0,i}\left(  \kappa\right)  r_{i}\\
&  =\sum_{i=0}^{\left\lfloor m/2\right\rfloor }r_{i}\left(  \left(
c_{i,0}+c_{i,1}\kappa\right)  g_{0,i}+\left(  c_{i,2}+c_{i,3}\kappa\right)
g_{0,i+1}\left(  \kappa\right)  +c_{i,4}g_{0,i+2}\left(  \kappa\right)
\right) \\
&  =\sum_{i=0}^{\left\lfloor m/2\right\rfloor +2}g_{0,i}\left(  r_{i}\left(
c_{i,0}+c_{i,1}\kappa\right)  +r_{i-1}\left(  c_{i-1,2}+c_{i-1,3}%
\kappa\right)  +r_{i-2}c_{i-2,4}\right)  .
\end{align*}
By symbolic computation
\[
\left(  c_{i,0}+c_{i,1}\kappa\right)  +\frac{r_{i-1}}{r_{i}}\left(
c_{i-1,2}+c_{i-1,3}\kappa\right)  +\frac{r_{i-2}}{r_{i}}c_{i-2,4}=0
\]
for $2\leq i\leq\left\lfloor \frac{m}{2}\right\rfloor $. Also
\begin{align*}
c_{0,0}+c_{0,1}\kappa &  =0,\\
r_{1}\left(  c_{1,0}+c_{1,1}\kappa\right)  +r_{0}\left(  c_{0,2}+c_{0,3}%
\kappa\right)   &  =0.
\end{align*}
The special cases at the top end of summation occur at $i=\frac{m}{2}$ or
$i=\frac{m-1}{2}$ (for $m$ being even or odd, respectively) and in fact
$c_{i,4}=0$ and $\left(  c_{i,2}+c_{i,3}\kappa\right)  +\frac{r_{i-1}}{r_{i}%
}c_{i-1,4}=0$ for these values of $i$.
\end{proof}

This completes the proof of Theorem \ref{big1}.

We use the same approach to Theorem \ref{big2}, where we set $i=m\geq0$ and
extract the common factor%
\begin{align*}
&  C_{m}=\frac{c\left(  2a_{1}+m-1,2a_{3}-m+1,2a_{3}-m+1,m+1\right)  \left(
-2a_{3}-1\right)  _{m}\left(  -2a_{2}-1\right)  _{m}}{s\left(  2a_{1}%
,2a_{3}+2,2a_{2}+2,0\right)  m!\left(  2a_{1}+1\right)  _{m}\left(
4\kappa+2n+1\right)  \left(  2\kappa+n\right)  }\\
&  =2^{2m+1}\frac{\left(  \kappa+a_{1}\right)  _{m}\left(  a_{1}+\frac{1}%
{2}\right)  _{m}\left(  -2a_{3}-1\right)  _{m}\left(  -2a_{2}-1\right)  _{m}%
}{\left(  -3-2\kappa-2a_{2}-2a_{3}\right)  _{2m+2}m!\left(  2a_{1}+1\right)
_{m}}.
\end{align*}
Next calculate (recall $n=a_{1}+a_{2}+a_{3}+1$)%
\begin{align*}
\frac{c\left(  2a_{1}+m+1,2a_{3}+1-m,2a_{2}+1-m,m+1\right)  }{c\left(
2a_{1}+m-1,2a_{3}+1-m,2a_{2}+1-m,m+1\right)  } &  =\frac{\left(  \kappa
+a_{1}+m\right)  \left(  2a_{1}+2m+1\right)  }{\left(  4\kappa+2n+1\right)
\left(  2\kappa+n\right)  },\\
\frac{c\left(  2a_{1}+m,2a_{3}+2-m,2a_{2}+2-m,m\right)  }{c\left(
2a_{1}+m-1,2a_{3}+1-m,2a_{2}+1-m,m+1\right)  } &  =\\
&  2\frac{\left(  \kappa+a_{2}+a_{3}+1-m\right)  \left(  \kappa+a_{2}%
+a_{3}+\frac{3}{2}-m\right)  }{\left(  4\kappa+2n+1\right)  \left(
2\kappa+n\right)  },\\
\frac{c\left(  2a_{1}+m+1,2a_{3}+1-m,2a_{2}+1-m,m-1\right)  }{c\left(
2a_{1}+m-1,2a_{3}+1-m,2a_{2}+1-m,m+1\right)  } &  =\frac{\left(
2\kappa+2a_{2}+1\right)  \left(  2\kappa+2a_{3}+1\right)  }{\left(
2a_{2}+1\right)  \left(  2a_{3}+1\right)  },\\
\frac{c\left(  2a_{1}+m,2a_{3}-m,2a_{2}+2-m,m\right)  }{c\left(
2a_{1}+m-1,2a_{3}+1-m,2a_{2}+1-m,m+1\right)  } &  =\frac{2\kappa+2a_{3}%
+1}{2a_{3}+1},\\
\frac{c\left(  2a_{1}+m,2a_{3}+2-m,2a_{2}-m,m\right)  }{c\left(
2a_{1}+m-1,2a_{3}+1-m,2a_{2}+1-m,m+1\right)  } &  =\frac{2\kappa+2a_{2}%
+1}{2a_{2}+1}.
\end{align*}
Write the right side of the summation formula in Theorem \ref{big2} as
$d_{m}=A_{m}F\left(  m;\kappa+1,a_{1},a_{2},a_{3}\right)  $. Then%
\begin{align*}
\frac{A_{m}}{C_{m}} &  =4\kappa\left(  \kappa+a_{1}+m\right)  \left(
4\kappa+3n+2\right)  \frac{\left(  2a_{2}-m+1\right)  \left(  2a_{3}%
-m+1\right)  }{\left(  2a_{2}+1\right)  \left(  2a_{3}+1\right)  },\\
\frac{A_{m-1}}{A_{m}} &  =\frac{2m\left(  2a_{1}+m\right)  \left(
\kappa+a_{2}+a_{3}+1-m\right)  \left(  \kappa+a_{2}+a_{3}+\frac{3}%
{2}-m\right)  }{\left(  \kappa+a_{1}+m\right)  \left(  2a_{1}+2m-1\right)
\left(  2a_{2}-m+1\right)  \left(  2a_{3}-m+1\right)  }.
\end{align*}
The expression $\left(  \frac{\left(  -2a_{2}-1\right)  _{m}\left(
-2a_{3}-1\right)  _{m}}{m!\left(  2a_{1}+1\right)  _{m}}t_{m}+d_{m-1}%
-d_{m}\right)  /C_{m}$ with $a_{i}$ replaced by $v_{i}$ (formal parameters
instead of integers) for $1\leq i\leq3$ becomes the left side of the equation
in the following.

\begin{theorem}
Suppose $v_{1},v_{2},v_{3}\notin-\frac{1}{2}+\mathbb{N}_{0}$, ($v_{0}%
=v_{1}+v_{2}+v_{3}+1$) and $m\in\mathbb{N}_{0}$ then
\begin{align*}
&  2\left(  2v_{1}+1+2m\right)  \left(  \kappa+v_{1}+m\right)  \left(
4\kappa+v_{0}\right)  F\left(  m+1;\kappa,v_{1},v_{2}+1,v_{3}+1\right) \\
&  -4\left(  4\kappa+v_{0}\right)  \left(  \kappa+v_{2}+v_{3}+\frac{3}%
{2}-m\right)  \left(  \kappa+v_{2}+v_{3}+1-m\right)  F\left(  m;\kappa
,v_{1},v_{2}+1,v_{3}+1\right) \\
&  -\left(  2v_{1}+m\right)  \left(  4\kappa+2v_{0}+1\right)  \left(
2\kappa+v_{0}\right)  F\left(  m+1;\kappa,v_{1}-1,v_{2}+1,v_{3}+1\right) \\
&  -\frac{\left(  2\kappa+2v_{2}+1\right)  \left(  2\kappa+2v_{3}+1\right)
}{\left(  2v_{2}+1\right)  \left(  2v_{3}+1\right)  }m\left(  4\kappa
+2v_{0}+1\right)  \left(  2\kappa+v_{0}\right)  F\left(  m-1;\kappa
,v_{1}+1,v_{2},v_{3}\right) \\
&  +\frac{\left(  2\kappa+2v_{3}+1\right)  }{\left(  2v_{3}+1\right)  }\left(
2v_{3}+1-m\right)  \left(  4\kappa+2v_{0}+1\right)  \left(  2\kappa
+v_{0}\right)  F\left(  m;\kappa,v_{1},v_{2}+1,v_{3}\right) \\
&  +\frac{\left(  2\kappa+2v_{2}+1\right)  }{\left(  2v_{2}+1\right)  }\left(
2v_{2}+1-m\right)  \left(  4\kappa+2v_{0}+1\right)  \left(  2\kappa
+v_{0}\right)  F\left(  m;\kappa,v_{1},v_{2},v_{3}+1\right) \\
&  +\frac{8m\kappa\left(  2v_{1}+m\right)  }{\left(  2v_{2}+1\right)  \left(
2v_{3}+1\right)  \left(  2v_{1}+2m-1\right)  }\left(  \kappa+v_{2}+v_{3}%
+\frac{3}{2}-m\right)  \left(  \kappa+v_{2}+v_{3}+1-m\right) \\
&  \times\left(  4\kappa+3v_{0}+2\right)  F\left(  m-1;\kappa+1,v_{1}%
,v_{2},v_{3}\right) \\
&  -\frac{4\kappa\left(  2v_{2}+1-m\right)  \left(  2v_{3}+1-m\right)
}{\left(  2v_{2}+1\right)  \left(  2v_{3}+1\right)  }\left(  4\kappa
+3v_{0}+2\right)  \left(  \kappa+v_{1}+m\right)  F\left(  m;\kappa
+1,v_{1},v_{2},v_{3}\right) \\
&  =0.
\end{align*}

\end{theorem}

\begin{proof}
The method is similar to the previous one; here we use $\gamma_{3}$ as point
of reference for the ratio calculations. Write the identity as $\sum_{n=1}%
^{8}\mu_{n}F\left(  \gamma_{n}\right)  =0$ and use the notation $t_{i}\left(
\gamma_{n}\right)  $ as before. As before $g_{i,k}\left(  \kappa\right)
:=\left(  \kappa+i\right)  _{k}\left(  -\kappa-v_{0}+i\right)  _{k}$ (but here
$v_{0}=\sum_{j=1}^{3}v_{j}+1$, different from the previous theorem). Set
\begin{align*}
p_{i}\left(  \kappa\right)   &  =\sum_{n=1}^{2}\mu_{n}\frac{t_{i+1}\left(
\gamma_{n}\right)  }{t_{i}\left(  \gamma_{3}\right)  }+\mu_{3}+\sum_{n=4}%
^{8}\mu_{n}\frac{t_{i}\left(  \gamma_{n}\right)  }{t_{i}\left(  \gamma
_{3}\right)  },i\geq1,\\
p_{0}\left(  \kappa\right)   &  =\sum_{n=1}^{2}\mu_{n}t_{1}\left(  \gamma
_{n}\right)  +\sum_{n=1}^{8}\mu_{n}.
\end{align*}
Then $p_{i}\left(  \kappa\right)  $ is a polynomial of degree 5 in $\kappa$.
By computer algebra we find the coefficients $\left\{  c_{ik}\right\}
_{k=0}^{5}$, which are functions of $i,m,v_{1},v_{2},v_{3}$, so that%
\[
p_{i}\left(  \kappa\right)  =c_{i,0}+c_{i,1}\kappa+\left(  c_{i,2}%
+c_{i,3}\kappa\right)  g_{i,1}\left(  \kappa\right)  +\left(  c_{i,4}%
+c_{i,5}\kappa\right)  g_{i,2}\left(  \kappa\right)  .
\]
Then the left side of the identity equals%
\begin{align*}
&  \sum_{i=0}^{\left\lfloor \left(  m+1\right)  /2\right\rfloor }p_{i}\left(
\kappa\right)  g_{0,i}\left(  \kappa\right)  r_{i}\\
&  =\sum_{i=0}^{\left\lfloor \left(  m+1\right)  /2\right\rfloor }r_{i}\left(
\left(  c_{i,0}+c_{i,1}\kappa\right)  g_{0,i}+\left(  c_{i,2}+c_{i,3}%
\kappa\right)  g_{0,i+1}\left(  \kappa\right)  +\left(  c_{i,4}+c_{i,5}%
\kappa\right)  g_{0,i+2}\left(  \kappa\right)  \right) \\
&  =\sum_{i=0}^{\left\lfloor \left(  m+1\right)  /2\right\rfloor +2}%
g_{0,i}\left(  r_{i}\left(  c_{i,0}+c_{i,1}\kappa\right)  +r_{i-1}\left(
c_{i-1,2}+c_{i-1,3}\kappa\right)  +r_{i-2}\left(  c_{i-2,4}+c_{i-2,5}%
\kappa\right)  \right)  ,
\end{align*}
where $r_{i}:=\dfrac{t_{i}\left(  \gamma_{3}\right)  }{g_{0,i}\left(
\kappa\right)  }$. The degree of $\kappa$ in $\mu_{1}$ and $\mu_{2}$ is 3 and
thus $c_{0,4}+c_{0,5}\kappa$ agrees with the generic value of $c_{i,4}%
+c_{i,5}\kappa$ at $i=0$. In the range $2\leq i\leq\left\lfloor \frac{m+1}%
{2}\right\rfloor $ the value of $\left(  c_{i,0}+c_{i,1}\kappa\right)
+\frac{r_{i-1}}{r_{i}}\left(  c_{i-1,2}+c_{i-1,3}\kappa\right)  +\frac
{r_{i-2}}{r_{i}}\left(  c_{i-2,4}+c_{i-2,5}\kappa\right)  $ is found to be
zero by symbolic computation. Also $c_{0,0}+c_{0,1}\kappa=0$ and $r_{1}\left(
c_{1,0}+c_{1,1}\kappa\right)  +r_{0}\left(  c_{0,2}+c_{0,3}\kappa\right)  =0$.
The special cases at the top end of summation occur at $i=\frac{m}{2}$ or
$i=\frac{m+1}{2}$ (for $m$ being even or odd, respectively) and in fact
$c_{i,4}+c_{i,5}\kappa=0$ and $\left(  c_{i,2}+c_{i,3}\kappa\right)
+\frac{r_{i-1}}{r_{i}}\left(  c_{i-1,4}+c_{i-1,5}\kappa\right)  =0$ for these
values of $i$.
\end{proof}

This completes the proof of Theorem \ref{big2}.

\section{Closing\ Comments}

Firstly we remark on the difficulty and complexity of the proof of the
intertwining property. We speculate that one reason is that the measure on the
linear transformations $\tau\left(  q\right)  $ is not uniquely defined. As
was seen in the presentation it is only the integrals of $P_{a,b}^{c}\left(
q\right)  $ that matter. By R\"{o}sler's result \cite{R1} for each
$x\in\mathbb{R}^{2}$ and $\kappa>0$ there is a positive Baire measure $\mu
_{x}$ such that $Vf\left(  x\right)  =\int fd\mu_{x}$ and the support of $\mu$
is contained in the closed convex hull of the $B_{2}$-orbit of $x$. The
explicit measure in Theorem \ref{Vint} requires $\kappa>\frac{3}{2}$ and is
not positive (combine the two parts into one and evaluate the kernel at
$u=\varepsilon,\phi_{1}=\phi_{2}=\frac{\pi}{2},\psi_{1}=\varepsilon,\psi
_{2}=\pi-\varepsilon$ for small positive $\varepsilon$; the value is $2\left(
2\kappa-3\right)  \varepsilon^{2}\cos\theta+\varepsilon^{4}\left(  1+\left(
7-\frac{16}{3}\right)  \cos\theta\right)  +O\left(  \varepsilon^{5}\right)  $;
let $\frac{\pi}{2}<\theta<\pi$). However the integral formula for the Bessel
function $K^{0}\left(  x,y\right)  $ is positive and works for $\kappa
>\frac{1}{2}$. It would be very interesting if one could develop explicit
results for the Coxeter groups $A_{n}$ or $B_{n}$ for $n\geq3$. The
corresponding compact Lie groups ($U\left(  n+1\right)  $ and $Sp\left(
n\right)  $, repectively) are fairly concrete but clearly we need more
powerful techniques than we used here for $n=2$.

Secondly there is the two-parameter problem. The group $B_{2}$ allows two
parameters in the associated differential-difference operators (see
(\ref{ddops})). The approach to finding the measure $\mu$ by starting with the
Lie group $Sp\left(  2\right)  $ did not (as yet) help us in the two-parameter
problem. This may be the hidden reason why the part of the intertwining
formula associated with the representations of types 1 and 2 (realized on
$x_{1}x_{2}$ and $x_{1}^{2}-x_{2}^{2}$ respectively) is so complicated
(involving derivatives); we were somehow (heuristically) close to the kernels
for $\kappa_{1}=\kappa+1,\kappa_{2}=\kappa$ and $\kappa_{1}=\kappa,\kappa
_{2}=\kappa+1$.

Finally we comment on the role of computer algebra. We produced $K_{n}\left(
x,y\right)  $ for $n\leq6$ by solving the equations $T_{i}^{x}K_{n}\left(
x,y\right)  =y_{i}K_{n-1}\left(  x,y\right)  ,i=1,2$. Then by use of the
symmetrized kernels $\frac{1}{8}\sum_{w\in B_{2}}K_{n}\left(  xw,y\right)  $
we conjectured that the measure $\mu$ is the right one. Borrowing the
polynomial $2\left(  q_{1}+q_{4}\right)  +q_{1}q_{4}-q_{2}q_{3}$ from the
$S_{3}$ paper \cite{D3} quickly helped to expand the conjecture. The part of
the intertwining operator dealing with types 1 and 2 (see Section \ref{sctB2})
turned out to be puzzling. Eventually we had a formula for $V$ that worked for
each $K_{n}\left(  x,y\right)  $ tried (up to $n=12$), but no proof. Then
experimenting with $P_{a,b}^{c}\left(  q\right)  $ led to discovering the
single-sum for $s\left(  \alpha\right)  $ (see Section \ref{sctsingl}). The
proof itself depends on classical transformations of hypergeometric series. In
turn the single sum expression (which had three free parameters) made feasible
enough experiments to formulate the identity in Theorem \ref{big1}. The
term-by-term calculations which prove the contiguity relations in Section
\ref{sctcont} would be very tedious without symbolic computation assistance,
and much space would be needed to write down every intermediate step. However
it is straightforward to verify these relations for any particular $m$
(magnitude depending on the size of the computer; $m=1\ldots6$ is not too
big). There is some satisfaction in coaxing a symbolic computation system
actually to prove a conjecture after it helped in the conjecture's formulation.


\begin{thebibliography}{99}                                                                                               %


\bibitem {B}W. Bailey, \emph{Generalized Hypergeometric Series}, Cambridge
University Press, Cambridge, 1935.

\bibitem {D1}C. Dunkl, Differential-difference operators associated to
reflection groups, \emph{Trans. Amer. Math. Soc.}\textbf{\ 311} (1989), 167-183.

\bibitem {D2}C. Dunkl, Operators commuting with Coxeter group actions on
polynomials, \emph{Invariant Theory and Tableaux}, (D. Stanton, ed.),
Springer, Berlin-Heidelberg-New York, 1990, pp. 107-117.

\bibitem {D3}C. Dunkl, Intertwining operators associated to the group $S_{3}$,
\emph{Trans. Amer. Math. Soc.}\textbf{\ 347} (1995), 3347-3374.

\bibitem {DJO}C. Dunkl, M. de Jeu, and E. Opdam, Singular polynomials for
finite reflection groups, \emph{Trans. Amer. Math. Soc.} \textbf{346} (1994), 237-256.

\bibitem {DX}C. Dunkl and Y. Xu, \emph{Orthogonal Polynomials of Several
Variables}, Encycl. of Math. and its Applications \textbf{81}, Cambridge
University Press, Cambridge, 2001.

\bibitem {H}S. Helgason, \emph{Groups and Geometric Analysis}, Academic Press,
New York, 1984.

\bibitem {R1}M. R\"{o}sler, Positivity of Dunkl's intertwining operator,
\emph{Duke Math. J.} \textbf{98} (1999), 445-463.

\bibitem {R2}M. R\"{o}sler, Dunkl operators: theory and applications,
\emph{Orthogonal Polynomials and Special Functions (Leuven 2002)}, (E. Koelink
and W. Van Assche, eds.), LNIM 1817, Springer, Berlin-Heidelberg-New York,
2003, pp. 93-135.

\bibitem {X}Y. Xu, A product formula for Jacobi polynomials, \emph{Special
Functions (Hong Kong 1999)}, (C. Dunkl, M. Ismail, R. Wong, eds.), World
Scientific, Singapore, 2000, pp. 423-430.
\end{thebibliography}
\end{document}